\documentclass[a4paper,12pt]{article}
\usepackage{cmap}
\usepackage[T1]{fontenc}
\usepackage[utf8]{inputenc}
\usepackage[a4paper]{geometry}
\usepackage{enumitem}
\usepackage{float}
\usepackage[unicode,hyperfootnotes=false,pdftex]{hyperref}
\usepackage{hyperxmp}
\usepackage{bookmark}
\usepackage{amsthm,amssymb}
\usepackage{mathtools}
\usepackage[cal=boondoxo]{mathalfa}
\usepackage{microtype}
\usepackage{xparse}
\usepackage{csquotes}
\usepackage{xcolor}
\usepackage[ISO]{diffcoeff}
\usepackage{booktabs}
\usepackage{array}
\usepackage[english]{babel}
\usepackage[useregional]{datetime2}
\usepackage{tikz}
\usepackage{graphicx}

\hypersetup{%
    pdftitle={Maximum likelihood estimator for skew Brownian motion: the convergence rate},
    pdfauthor={Antoine Lejay; Sara Mazzonetto},
    pdfdate={2023-01-02},
    pdfkeywords={Skew Brownian motion; maximum likelihood estimator (MLE); statistical estimation; null recurrent process.},
    pdflang={en},
    pdfmetalang={en},
    pdfsubject = {Inference of the Skew Brownian motion},
}

\SetEnumitemKey{thm}{nolistsep,topsep=-\parskip,label={\textnormal{(\roman*)}}}
\SetEnumitemKey{compact}{noitemsep,topsep=-\parskip,leftmargin=0em,labelindent=2em}

\newcommand{\cvstable}[1][]{\xrightarrow[#1]{\text{stably}}}
\newcommand{\cvproba}[1][]{\xrightarrow[#1]{\text{prob.}}}

\newcommand{\RR}{\mathbb{R}}
\newcommand{\EE}{\mathbb{E}}
\newcommand{\PP}{\mathbb{P}}
\newcommand{\NN}{\mathbb{N}}

\newcommand{\cN}{\mathcal{N}}
\newcommand{\cF}{\mathcal{F}}
\newcommand{\cL}{\mathcal{L}}
\newcommand{\cs}{\mathcal{s}}

\newcommand{\vd}{\,\mathrm{d}}

\newcommand{\sP}{\mathsf{P}}
\newcommand{\sd}{\mathsf{d}}
\newcommand{\sa}{\mathsf{a}}
\newcommand{\sD}{\mathsf{D}}

\newcommand{\uPsi}{\underline{\Psi}}

\newcommand{\ind}[1]{\mathbf{1}_{#1}}

\newcommand{\zig}{\tikz \draw[thick] (0,0) -- (3pt,3pt) -- (6pt,0pt) -- (9pt,3pt);}

\newcommand\given{\nonscript\:\delimsize\vert\nonscript\:\mathopen{}} 
\newcommand\SetSymbol[1][]{\nonscript\:#1\vert\nonscript\:\mathopen{}\allowbreak}
\DeclarePairedDelimiterX\Set[1]\{\}{%
  \renewcommand\given{\SetSymbol[\delimsize]}#1}
\DeclarePairedDelimiterX\Prb[1]{[}{]}{%
  \renewcommand\given{\SetSymbol[\delimsize]}#1}
\DeclarePairedDelimiterX\Paren[1](){#1}
\newcommand{\coef}[1]{\mathsf{#1}}
\DeclarePairedDelimiter{\abs}{|}{|}
\DeclarePairedDelimiter{\Floor}{\lfloor}{\rfloor}

\DeclareMathOperator{\Cov}{Cov}
\DeclareMathOperator{\atanh}{atanh}
\DeclareMathOperator{\erfc}{erfc}
\DeclareMathOperator{\erfcx}{erfcx}

\DeclareMathOperator{\sgn}{sgn}
\DeclareMathOperator*{\argmax}{arg\,max}

\newcommand{\eqlaw}{\stackrel{\mathrm{law}}{=}}

\newtheorem{theorem}{Theorem}
\newtheorem{proposition}{Proposition}
\newtheorem{lemma}{Lemma}
\newtheorem{corollary}{Corollary}

\theoremstyle{definition}
\newtheorem{definition}{Definition}
\newtheorem{notation}{Notation}
\newtheorem{data}{Data}

\theoremstyle{remark}

\newtheorem{remark}{Remark}

\begin{document}

\title{Maximum likelihood estimator\\
for skew Brownian motion:
\\
the convergence rate}
%{Estimation of the rate of convergence of the MLE estimator for the SBM}

\author{Antoine Lejay\footnote{Universit\'e de Lorraine, CNRS, IECL, Inria, F-54000 Nancy, France;
  ORCID: \href{https://orcid.org/0000-0003-0406-9550}{0000-0003-0406-9550}; \texttt{Antoine.Lejay@univ-lorraine.fr}}
    \and Sara Mazzonetto\footnote{Universit\'e de Lorraine, CNRS, IECL, Inria, F-54000 Nancy, France; 
	ORCID: \href{https://orcid.org/0000-0001-6187-2716}{0000-0001-6187-2716};
    \texttt{Sara.Mazzonetto@univ-lorraine.fr}}
}

\date{\today}

\maketitle

\begin{abstract}
We give a thorough description of the asymptotic property of the
maximum likelihood estimator (MLE) of the skewness parameter of a 
Skew Brownian Motion (SBM). Thanks to recent results on the Central
Limit Theorem of the rate of convergence of estimators for the SBM,
we prove a conjecture left open that the MLE has asymptotically a mixed normal
distribution involving the local time with a rate of convergence of order $1/4$. 
We also give a series expansion of the MLE and study the asymptotic behavior of the
score and its derivatives, as well as their variation with the skewness parameter. 
In particular, we exhibit a specific behavior when the SBM is actually 
a Brownian motion, and quantify the explosion of the coefficients of the expansion when the skewness
parameter is close to $‐1$ or $1$. 
\end{abstract}

\begin{quote}

\end{quote}

%\tableofcontents

% Intro: TO BE WRITTEN. Here, possible structure: \\
% SBM introduction with simple description and motivation based on importance in applications. \\
% Recall literature on estimation of SBM and cite new result which allow this work. \\
% Say what we do here: expansion, rates, constants, study of the coefficients and the fact that they are suitable for MC. \\
% Applications and numerics: hypothesis testing and numerics for skewness, coefficients, constants.\\
% Mention the fact that for parameters close to $\pm 1$ there is a problem in the expansion and also in the convergence of the pivotal quantity but mention also the remarks in~\cite{lejay2019} about it.\\
% Conclude with short outline.

% \bigskip

Skew and other singular diffusions attract more and more interest in 
modeling diffusive stochastic behavior in presence of semi-permeable barriers, 
discontinuities, and thresholds. Beyond theoretical studies, simulation and 
inference are also necessary tools for practical purposes.
For some example of applications in various fields, see 
\textit{e.g.} \cite{mota14a,decamps,ramirez2,zhang,ovaskainen03a,cantrell99a,gairat17a,bounebache_2014,ILM_22}
among others.

The inference of skew diffusion cannot follow from a simple adaptation
of known techniques for Stochastic Differential Equations (SDE)
as the ones presented in~\cite{ih,K}.
In fact their distributions are singular with the ones of classical SDE. 
The limits are usually mixed normal ones and the rate is not necessarily $1/2$.
The work dealing with the inference of skew diffusion is rather limited. 
Let us cite however \cite{bardou,lejay2014,lejay2019,mazzonetto19a,pigato,lejay_pigato2020,lejay_pigato2019,lejay_pigato2018,su_chan2015,su_chan2017} for frequentist inference 
and \cite{barahona,araya,Bai2021} for Bayesian inference.

The Skew Brownian motion (SBM) is a basic brick for constructing 
Skew diffusion, as several transformations reduces Skew diffusions
to SBM~\cite{lejay_sbm}. This latter process depends on a single parameter $\theta\in[-1,1]$
---~the \emph{skewness parameter}~--- which rules out its behavior
when it crosses zero. For $\theta=0$, the SBM is a Brownian motion.
For $\theta=\pm 1$, it is a Reflected Brownian motion.

A series of works considers the inference of the skewness parameter 
from high-frequency observations. In~\cite{lejay2014}, the authors have 
given an asymptotic expansion of the Maximum Likelihood Estimator (MLE)
$\theta_n$ of $\theta$ around $\theta=0$ in power of $n^{-1/4}$, where $n$ is the number
of observations. A heuristic explanation of the power $1/4$ is given by analogy
with the Skew Random Walk~\cite{lejay2018}, where
the MLE depends on the local time at zero of the discrete walk, 
which a random quantity of order~$\sqrt{n}$. 
Indeed in~\cite{lejay2019}, where the consistency of the MLE and another estimator is proved,
it was empirically observed that the rate of convergence should be $1/4$, 
meaning that the observed points that \textquote{carry the information}
are those close to~$0$, and are of order~$\sqrt{n}$.
This also explains 
why the asymptotic limit of the MLE involves the local time at~$0$.
The limit is  
of type $\cs(\theta)G/\sqrt{L_1}$ (when the process is observed on $[0,1]$),
where $(L_t)_{t\geq 0}$ is the symmetric local time at $0$ of the SBM and 
$G$ is a centered unit, Gaussian independent from the SBM.
The value of~$\cs(\theta)$ was empirically observed as closed to $\kappa\sqrt{1-\theta^2}$.
It was conjectured in~\cite{lejay2019} that~$\cs(\theta)=\kappa\sqrt{1-\theta^2}$
as for the Skew Random Walk. 

The article~\cite{mazzonetto19a} brought the missing result required to establish
a Central Limit Theorem on the MLE and other estimators. 
A related result may also be found in~\cite{robert_2022}. As
for the results in~\cite{lejay2019}, this is based on an extension 
of the work of J.~Jacod~\cite{jacod98} that cannot be applied directly 
as the SBM has a singular distribution with respect to the one of the BM.

In this article, we first prove the asymptotic normality of the MLE
of the skewness parameter $\theta$ of the SBM.
More precisely, we establish that 
\begin{equation*}
    n^{1/4}(\theta_n-\theta)\cvstable[n\to\infty]{} \cs(\theta)W(L_T)
    \eqlaw \cs(\theta)\frac{W(1)}{\sqrt{L_T}}
\end{equation*}
for a Brownian motion $W$ independent from the SBM. 
In particular, we show that~$\cs(\theta)=\kappa(\theta)\sqrt{1-\theta^2}$ 
with a pre-factor~$\kappa(\theta)$ that varies slowly.

Second, using a recent asymptotic inversion formula~\cite{lm22}, we 
establish a series expansion of $\theta_n$ as 
\begin{equation}
    \label{eq:intro:1}
    \theta_n
   =\theta+\sum_{k=1}^{+\infty} \sD_{k,n}(\theta) 
   \Paren*{
   \frac{  \cs(\theta) \sP_n(\theta)}{T^{1/4} n^{1/4}}}^{k},
\end{equation}
for some random quantities $\sD_{k,n}(\theta)$ whose
asymptotic behavior is also studied, and~$\sP_n(\theta)$
is asymptotically pivotal (\textit{i.e.},~the asymptotic law does not depend on the parameter).
All these
results are based on the asymptotic behavior
of the score and their derivatives at any orders. 
Expansion~\eqref{eq:intro:1} provides some insight on the behavior 
of the MLE in function of the number of samples $n$ 
and the true value of~$\theta$. In particular, the closer~$\theta$ is to~$\pm 1$, 
the more skewness is observed. 

Furthermore, we show that $\sD_{2k+1,n}(\theta=0)$ converges
to $0$ so that a more precise expression of the involving a multivariate mixed Gaussian distribution
could be given for the MLE of~$\theta=0$. This expression is an alternative to the
one already given in~\cite{lejay2014}. 

Moreover, we give some numerical experiments on the rate of convergence
of the distribution functions of the MLE and the score and its derivatives 
towards their limiting distributions. The rate of convergence of this 
Berry-Esseen type analysis seems to depend on~$\theta$. 
This will be subject to 
further study.

Finally, we study the behavior of the limiting coefficients in function of 
$\theta$. In particular, the expansion in \eqref{eq:intro:1} exhibits
a boundary layer estimate for $\theta$ close to~$\pm 1$ as the coefficients
explode in powers of~$(1-\theta^2)^{-1/2}$.

\bigskip

\textbf{Outline.} We present our main results in Section~\ref{sec:main}. 
They are built on the asymptotic behavior of the score and its derivative
which we give in Section~\ref{sec:score}. 
We also study numerically the rate of convergence of the scores
and their derivatives towards their limits in Section~\ref{sec:rate}.
The properties of 
the limiting coefficients are studied in Section~\ref{sec:score_coefficients}. 
And the proofs of our main theorems are given in Section~\ref{sec:proofs}.

%%%%%%%%%%%%%%%%%%%%%%%%%%%%%%%%%%%%%%%%%%%%%%%%%%%%%%%%%%%%%%%%%%%%%%
%%%%%%%%%%%%%%%%%%%%%%%%%%%%%%%%%%%%%%%%%%%%%%%%%%%%%%%%%%%%%%%%%%%%%%
\section{Main results}

\label{sec:main}

The SBM $X$ of parameter $\theta\in[-1,1]$ solves the SDE 
\begin{equation}
    \label{eq:sde}
X_t= B_t+ \theta L_t,\ t\geq 0
\end{equation}
for $B$ a Brownian motion and $\Set{L_t}_{t\geq 0}$ its symmetric local time of $X$
at point~$0$.
Actually, the SDE \eqref{eq:sde} has a unique strong solution~\cite{lejay_sbm,harrison}.
No solution exists when~$\abs{\theta}>1$. For $\theta=\pm 1$, the SBM is a (positively 
if $\theta=1$ or negatively if $\theta=-1$) Reflected Brownian motion.

We denote by $(\Omega,\cF,\PP_\theta)$ the underlying probability space
and by $\Set{\cF_t}_{t\geq 0}$ the filtration with respect
to which $X$ is adapted. 
This filtration may be taken as the
one generated by the driving Brownian motion 
and may be assumed to satisfy the usual conditions (i.e.~it is complete and right continuous).

We are concerned with the estimation of the parameter of the SBM observed
at discrete times over a finite time window. We will establish limits
in high-frequency.

\begin{data}
    \label{data:1}
    We observed the SBM $\Set{X_t}_{t\in[0,T]}$ of parameter $\theta_0\in(-1,1)$ 
    at times $\Set{t_i}_{i=0,\dotsc,\Floor{nT}}$ on a time window $[0,T]$ with $t_i:=i/n$. 
    The starting point is $X_0=0$. Note that $n$ is the number of sample per unit of time.
\end{data}
\begin{remark}
    We impose $X_0=0$ to ensure that $L_T\neq 0$. When $X_0\neq 0$, 
    we could still apply our methodology on a random window $[\tau_0\wedge T,T]$
    where $\tau_0$ is the first hitting time from $0$. If $\tau_0>T$, 
    then no observation can be used to estimate $\theta$.
\end{remark}

The density transition function of the SBM of parameter $\theta$ is \cite{walsh,lejay_sbm}
\begin{equation*}
    p_\theta(t,x,y) := 
    \frac{1}{\sqrt{2\pi t}}
    \begin{cases}
	\exp\Paren*{-\frac{(x-y)^2}{2t}}+\theta \exp\Paren*{-\frac{(x+y)^2}{2t}}&\text{ if }x\geq 0,\ y\geq 0,\\
	\exp\Paren*{-\frac{(x-y)^2}{2t}}-\theta \exp\Paren*{-\frac{(x+y)^2}{2t}}&\text{ if }x\leq 0,\ y\leq 0,\\
	(1-\theta)\exp\Paren*{-\frac{(x-y)^2}{2t}}&\text{ if }x\geq 0,\ y\leq 0,\\
	(1+\theta)\exp\Paren*{-\frac{(x-y)^2}{2t}}&\text{ if }x\leq 0,\ y\geq 0.
    \end{cases}
\end{equation*}
The SBM is null recurrent process with invariant measure $\mu_\theta( \vd x) :=\mu_\theta(x) \vd x$ with
\begin{equation*}
    \mu_\theta(x)
    :=\begin{cases}
	1+\theta&\text{ if }x\geq 0,\\
	1-\theta&\text{ if }x<0.
    \end{cases}
\end{equation*}

When observed at regular times as in Data~\ref{data:1}, we call \emph{likelihood} the random function:
\begin{equation*}
    \Lambda_n(\theta)=\prod_{i=0}^{\Floor{nT}-1} p_\theta(\Delta t,X_{t_i},X_{t_{i+1}}),\ \theta\in[-1,1]
    \text{ with }\Delta t=\frac{1}{n}.
\end{equation*}
Since $[-1,1] \ni \theta \mapsto p_\theta(\Delta t,x,y)$ is analytic, $\theta\mapsto \Lambda_n(\theta)$ is also analytic.
The \emph{score} is %$S(n,\theta):=$
$\partial_\theta \log \Lambda_n(\theta)$. 

Let us define
\footnote{In the definition of $k_\theta(x,y)$, there is a mistake in \cite{lejay2019}: 
$\Delta t$ has to be replaced by $T$. Note that in \cite{lejay2019}, $n$ is the number of samples
while here, $n$ is the number of samples per unit. Therefore, the limit distribution 
in \cite{lejay2019} is the local time at time $1$ of $\sqrt{T}X$, while here it is the one of $X$.}
\begin{equation*}
    k_\theta(x,y)
    :=\partial_\theta \log p_\theta(1,x,y)
    =
    \frac{\partial_\theta p_\theta(1,x,y)}{p_\theta(1,x,y)}
    =\frac{\sgn(y)}{\sgn(y)\theta+\exp(2(xy)^+)},
\end{equation*}
where $(xy)^+$ stands for $\max(0, xy)$.
For any $n>0$, the following scaling holds true:
\begin{equation*}
    k_\theta(x\sqrt{n},y\sqrt{n})
    =\partial_\theta \log p_\theta(\Delta t,x,y)\text{ with }\Delta t=\frac{1}{n}.
\end{equation*}

\begin{remark} \label{rem:score:k}
The score rewrites %$S(n,\theta)$ 
$\partial_\theta \log \Lambda_n(\theta) = \sum_{i=0}^{\Floor{nT}-1} k_\theta(X_{t_i}\sqrt n,X_{t_{i+1}}\sqrt n)$.
\end{remark}

\begin{proposition}[{{\cite{lejay2019}}}]
    The maximum likelihood estimator (MLE) 
	\[\theta_n:=\argmax_{\theta\in[-1,1]} \Lambda_n(\theta)\]
    is the unique solution to $\partial_\theta \log \Lambda_n(\theta_n)=0$ %$S(n,\theta_n)=0$ 
    and is a consistent estimator of the parameter $\theta_0$ of the SBM under $\PP_{\theta_0}$. 
\end{proposition}

The goal of this paper is to refine the latter result providing asymptotic information.

\begin{figure}
    \begin{center}
    \includegraphics{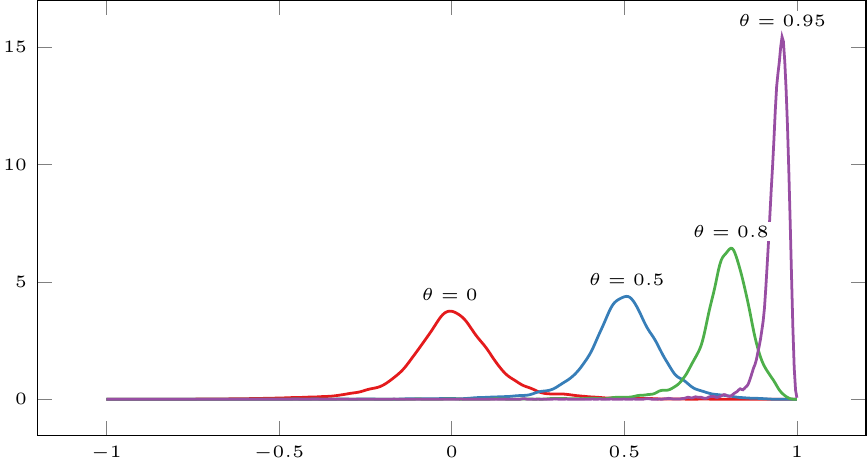}
    \caption{\label{fig:density:mle} Density of the MLE for various values of $\theta$ using $10\,000$
    samples of the SBM using $\Delta t=10^{-4}$ for $T=1$.}
\end{center}
\end{figure}

In Figure~\ref{fig:density:mle}, we plot the empirical density of the MLE for various
values of $\theta$. We use the method in \cite{lejay12} for the simulation of the SBM, 
while the MLE associated to each trajectory is obtained maximizing numerically the log-likelihood.
We see that the MLE is concentrated around the true value. 
The more $\theta$ is closer to $1$, the more the density is skewed to 
the left.
In the next section we examine this behavior.
For reasons related to symmetry, for the remainder of the document, we restrict to consider the case $\theta\geq 0$.

%%%%%%%%%%%%%%%%%%%%%%%%%%%%%%%%%%%%%%%%%%%%%%%%%%%%%%%%%%%%%%%%%%%%%%
%%%%%%%%%%%%%%%%%%%%%%%%%%%%%%%%%%%%%%%%%%%%%%%%%%%%%%%%%%%%%%%%%%%%%%
\subsection{Asymptotic mixed normality of the MLE estimator}

We establish, in the forthcoming Theorem~\ref{thm:mixed_normality}, the asymptotic mixed normality of the MLE, 
with the non standard rate of $n^{1/4}$. 
This result has already been proven, when $\theta=0$, in~\cite{lejay2014}. 
Before to state the result, 
let us first introduce a definition specifying the nature of the asymptotic normality.

\begin{definition}[Class of $\cL$-mixed normal distribution]
    Let $L_1$ be the local time at point $0$ and time $1$ of the SBM $X$ and $G\sim\cN(0,1)$
    be independent from $\cF_1$ (hence from~$L_1$). Note that 
		the distribution of $L_1$ does not depend on $\theta$.
    A random variable $M$ is said to be \emph{$\cL$-mixed normal distributed} if 
    \begin{equation*}
	M\eqlaw \frac{G}{\sqrt{L}_1}.
    \end{equation*}
\end{definition}

Note that a $\cL$-mixed normal distribution has an infinite second moment.
It is a symmetric, unimodal distribution with heavy tails.

\begin{remark}[Simulation of the local time]
The local time of the SBM is equal in distribution to the one of the Brownian motion, hence $L_1\eqlaw\abs{H}$ for $H\sim\cN(0,1)$.
Therefore, the local time and $\cL$-mixed normal distributions are easily simulated.
%
%Alternatively, as the starting point is $0$, the Paul Lévy theorem and the results in \cite{lepingle95a}
%implies that
%\begin{equation*}
%    L_1\eqlaw \frac{1}{2}\Paren*{U+\sqrt{V+U^2}}
%    \text{ with }U\sim\cN(0,1)\text{ and }V\sim\cE(1/2),
%\end{equation*}
%where $U$ and $V$ are independent.
%
\end{remark}
\begin{remark}[Scaling]
    \label{rem:scaling}
    Let $L$ be the symmetric local time of a SBM $X$ at the point $0$ and $W$ be
     a Brownian motion $W$ independent from~$X$. 
	They both satisfy a scaling property, in particular for the local time it holds $L_T\eqlaw \sqrt{T}L_1$.
     Hence, 
    \begin{equation*}
	\frac{W(L_T)}{L_T} \eqlaw \frac{W(1)}{\sqrt{L_T}}
	\eqlaw \frac{W(1)}{T^{1/4}\sqrt{L_1}}.
    \end{equation*}
    In other words, $T^{1/4}W(L_T)/L_T$ and $T^{1/4}G/\sqrt{L_1}$
    (for $G\sim\cN(0,1)$ independent from~$X$)
    are $\cL$-mixed normal distributed.
\end{remark}

We also introduce a quantity related to the asymptotic variance of the estimator and 
to the Fisher information:
\begin{equation} \label{eq:s_theta:def}
	    \cs(\theta):= \Paren*{ - \iint_{\RR^2} \mu_\theta(x) p_\theta(1,x,y)  \partial_\theta k(x,y) \vd x \vd y}^{\!-1/2} \in \RR.
\end{equation}
The quantity $\cs(\theta)$ will be studied in more details in Proposition~\ref{prop:s_theta}, 
Remark~\ref{rem:cs},
and in Section~\ref{sec:score_coefficients}, in particular in Section~\ref{sec:amn}.

We are now ready to state the main result which says that
$(nT)^{1/4}\cs(\theta)^{-1}(\theta_n-\theta)$
is asymptotically a $\cL$-mixed normal distribution under $\PP_{\theta}$.

\begin{theorem}[Asymptotic $\cL$-mixed normality of the MLE for the SBM]
    \label{thm:mixed_normality}
    Let $\theta\in (-1,1)$. The MLE estimator $\theta_n$ is asymptotically mixed normal
    with 
    \begin{equation*}
	n^{1/4}(\theta_n-\theta)\cvstable[n\to\infty] \cs(\theta)\frac{H}{\sqrt{L_T}},
	\text{ under }\PP_{\theta},
    \end{equation*}
    where $H\sim\cN(0,1)$ independent from $\cF_T$ (hence of $L_T$).
\end{theorem}

The proof of Theorem~\ref{thm:mixed_normality} 
follows from the results of the next section,
Section~\ref{sec:asym},
which rely on the study of the score and its derivatives that we propose in Section~\ref{sec:score}.

Heuristically, the non standard rate of $n^{1/4}$, can be explained
by the fact that the quality of the estimation depends mainly 
on the time spent by the SBM around $0$, and that the fraction
of observations when $\Set{X_{t_i}}_{i=0,\dotsc,n}$ is of order~$\sqrt{n}$.
This fact is rigorously established for the Skew Random Walk where 
the local time is really the occupation time at the point where the bias-dynamic 
is perturbed~\cite{lejay18}.

Although it was conjectured in \cite{lejay2019} from the results on the Skew Random Walk
that the coefficient $\cs(\theta)$ in front of the mixed Gaussian 
should be proportional to~$\sqrt{1-\theta^2}$, 
we find a slowly varying pre-factor.% which will be studied extensively in Section~\ref{sec:amn}.

\begin{proposition} 
    \label{prop:s_theta}
For all $\theta \in (-1,1)$, the function $(x,y) \mapsto \mu_\theta(x) p_\theta(1,x,y)  \partial_\theta k(x,y) $ is integrable 
and its integral is negative, so $\cs(\theta)\in (0,\infty)$ is well defined.
Besides there exist two real constants $0< c_1 \leq c_2 <\infty$ such that for all $\theta\in(-1,1)$
	\begin{equation*}
		c_-\sqrt{1-\theta^2}\leq \cs(\theta)\leq c_+\sqrt{1-\theta^2}.
    \end{equation*}
\end{proposition}

We find that $0.79 \leq c_- < c_+\leq 0.88$, 
which is consistent with the numerical observations of \cite{lejay2019}.

The proof of Proposition~\ref{prop:s_theta} is provided in Section~\ref{sec:amn}
where a more precise statement is formulated.
Actually we show in Remark~\ref{rem:stheta} that an accurate approximation of~$\cs(\theta)$ is given by
\begin{equation*}
    \cs(\theta)  
    \approx \frac{\sqrt{1-\theta^2}}{
\sqrt{1.292+ 0.232 \, \theta^2+  0.071 \, \theta^4}}.
\end{equation*}

%%%%%%%%%%%%%%%%%%%%%%%%%%%%%%%%%%%%%%%%%%%%%%%%%%%%%%%%%%%%%%%%%%%%%%
%%%%%%%%%%%%%%%%%%%%%%%%%%%%%%%%%%%%%%%%%%%%%%%%%%%%%%%%%%%%%%%%%%%%%%
\subsection{Asymptotic expansion for the MLE estimator} 
\label{sec:asym}

Let us first consider the following family of statistics of interest: 
\begin{align}
     \label{eq:coefS}
%     \coef{S}_{0}(n,\theta):=\frac{1}{n^{1/4}}\sum_{i=0}^{n-1}  k_\theta(X_{t_i}\sqrt{n},X_{t_{i+1}}\sqrt{n}),\\
     \coef{S}_{m}(n,\theta):=\frac{1}{n^{1/2}}\sum_{i=0}^{\Floor{nT}-1} \partial_\theta^{m} k_\theta(X_{t_i}\sqrt{n},X_{t_{i+1}}\sqrt{n})\text{ for }m\geq 0.
 \end{align}
Remark~\ref{rem:score:k} shows that
$\coef{S}_{m}(n,\theta) = \partial_\theta^m \Paren*{\partial_\theta \log \Lambda_n(\theta)}/n^{1/2}$
so that $\coef{S}_{m}(n,/\theta)$ is the renormalized $m$-th order derivative of the score.

Let us also consider the statistics, on $\{\coef{S}_1(n,\theta) \neq 0\}$,
\begin{align*}
    \sd_{k,n}(\theta)  := 
    \frac{-1}{k!}\frac{\coef{S}_k(n,\theta)}{\coef{S}_1(n,\theta)}
    \text{ for }k\geq 0
	\quad \text{ and } \quad 
	\sP_n(\theta):= n^{1/4} \frac{T^{1/4}}{\cs(\theta)}  \sd_{0,n}(\theta),
\end{align*}
as well as the constants
\begin{equation}
	\label{eq:xi:def}
	\xi_m(\theta):=
	\iint_{\RR^2} \mu_\theta(x) p_\theta(1,x,y)\partial_\theta^{m}k(x,y)\vd x \vd y, \quad \text{ for } m\in \NN.
	%\\=
	%m!(-1)^m \iint_{\RR^2} \mu_\theta(\vd x)p_\theta(T,x,y)k^{m+1}(x,y)\vd x \vd y.
\end{equation}
The integrability of $(x,y) \mapsto \mu_\theta(x) p_\theta(1,x,y)\partial_\theta^{m}k(x,y)$ is shown in Section~\ref{sec:score_coefficients}, together with other properties of $\xi_m(\theta)$.

\begin{remark}
    \label{rem:cs}
    The quantity $\cs(\theta)$ 
    given by \eqref{eq:s_theta:def}
    is related to $\xi_1(\theta)$: $\cs(\theta)^{-2} = -\xi_1(\theta)$.
\end{remark}
%%

%%%
\begin{proposition}
\label{prop:pn_dn}
Under $\PP_\theta$,
\begin{enumerate}[thm]
    \item the statistics $\sP_n(\theta)=n^{1/4}T^{1/4}\sd_{0,n}(\theta)/\cs(\theta)$
	is asymptotically $\cL$-mixed normal distributed;
    \item for every $k\geq 1$,
\begin{equation} 
    \label{eq:def:dn}
    \sd_{k,n}(\theta)\cvproba[n\to\infty] \sd_k(\theta):=\frac{-\xi_k(\theta)}{k!\xi_1(\theta)}
    = \frac{\cs(\theta)^2}{k!}\xi_k(\theta). 
\end{equation}
\end{enumerate}
Moreover, under $\PP_0$, for every $k\geq 1$, $\sd_{2k,n}(0)\cvproba[n\to\infty] 0$ and
$\sd_{2k}(0)$
is, up to a multiplicative constant, $\cL$-mixed normal distributed.
Furthermore  $\Set{n^{1/4} \sd_{2k,n}(0)}_{k=0,\dotsc,m}$ converges stably
for any $m\geq 1$. The limit is identified in Proposition~\ref{prop:score} in Section~\ref{sec:score}.
\end{proposition}

\begin{proof}
    This is an immediate consequence of Proposition~\ref{prop:score} in Section~\ref{sec:score}
    on the asymptotic behavior of the score combined
    with Remark~\ref{rem:cs}.
\end{proof}

In Theorem~\ref{thm:asym}, we give an expansion of the MLE in term of $\sd_{0,n}(\theta)$ (and so of $\frac{\cs(\theta) \sP_n(\theta)}{T^{1/4} n^{1/4}}$). 
It follows from applying Theorem~3 in~\cite{lm22}.
For $\theta=0$, 
such a type of expansion was already given in~\cite{lejay2014} in the form provided in equation~\eqref{eq:mle:neq0} below. We also provide an alternative expansion based on a finer analysis of the coefficients' asymptotic behavior.

Let us introduce some notation:
for any $m\in \NN\cup \{\infty\}$,
let 
\begin{equation}
	\Phi^{[m]}(\delta,x):= \sum_{k=1}^m \delta_k x^k
\end{equation}
be the formal power series in $x$ of coefficients $\delta=\Set{\delta_k}_{k\geq 1}$.

\begin{theorem}[{Asymptotic expansion, cf.~\cite[Theorem~3]{lm22}}]
    \label{thm:asym}
    For any integer $n\geq 1$, any $N\geq 1$, and any $\theta \in (-1,1)$, the MLE satisfies under $\PP_{\theta}$,  
\begin{equation}
    \label{eq:mle:neq0}
   \theta_n
   =
	\theta + \Phi^{[\infty]}( \sD_{\cdot,n}(\theta), \sd_{0,n}(\theta))
   =
	\theta+\sum_{k=1}^{+\infty} \sD_{k,n}(\theta) \sd_{0,n}(\theta)^k,
   %=
   %\theta+\sum_{k=1}^{+\infty} \sD_{k,n}(\theta) 
   %\Paren*{
   %\frac{  \cs(\theta) \sP_n(\theta)}{T^{1/4} n^{1/4}}}^{k},
   \end{equation}
 where $\sD_{1,n}(\theta):=1$ and $\sD_{q,n}(\theta)$, $q =2,\ldots, N$, are given by the following recursive formula:
 \begin{equation} 
     \label{eq:sD:recursive}
	\sD_{q,n}(\theta) :=  \sum_{m=2}^{q} \sd_{m,n}(\theta) \sum_{k_1+\dotsb+k_m=q}
    \sD_{k_1,n}(\theta)\dotsb \sD_{k_m,n}(\theta).
\end{equation}
If $\theta=0$, another expansion holds
\begin{equation}
    \label{eq:mle:0}
    \theta_n 
	= \Phi^{[\infty]}(\sa_{\cdot,n}, n^{-1/4})
	= \sd_{0,n}(0) + \sum_{\substack{ k\geq 3\\ k \text{ odd}}} \frac{\sa_{k,n}}{n^{k/4}},
\end{equation}
where $\sa_{k,n}$ for $k =0,\ldots, N$ is $\sa_{2 q,n}:=0$ if $q\geq 0$,  and
\begin{gather*}
    \sa_{1,n}:=n^{1/4}\sd_{0,n}(0) = \cs(0) \sP_n(0)/{T^{1/4}},
    \\
\begin{multlined}  
	\sa_{2q+1,n}  
	  := \sum_{\substack{ m=2 \\ m \text{ even }}}^{2q}  
	 n^{1/4} \sd_{m,n}(0) \sum_{k_1+\dotsb+k_m=2q} \sa_{k_1,n} \cdots \sa_{k_m,n} 
	\\	
	  + \sum_{\substack{ m=3 \\ m \text{ odd }}}^{2q+1} \sd_{m,n}(0)
	 \sum_{k_1+\dotsb+k_m = 2q+1} \sa_{k_1,n} \cdots \sa_{k_{m},n}.
\end{multlined}
\end{gather*}
\end{theorem}
\begin{remark}
From~\cite{lm22}, it can be seen that $\Phi(\sd_{\cdot,n}(\theta),\Phi(\sD_{\cdot,n}(\theta),x))=-x$. 
For the first coefficients, it holds 
    \begin{gather*}
	\sD_{2,n}=\sd_{2,n},\ 
	\sD_{3,n}=\sd_{3,n} + 2 \sd_{2,n}^2,
	\\
	\sD_{4,n}=\sd_{4,n} + 5 \sd_{2,n}\sd_{3,n} + 5 \sd_{2,n}^3
	\\
	\text{and }
	\sD_{5,n}= \sd_{5,n} + 6 \sd_{2,n}\sd_{4,n} + 21 \sd_{2,n}^2\sd_{3,n} + 14\sd_{2,n}^4 +3 \sd_{3,n}^2.
    \end{gather*}
In the case of the sequence in~\eqref{eq:mle:0},
 	\begin{gather*}
	\sa_{3,n}=\sd_{3,n}(0) \sP_n(0)^3 \cs(0)^3/T^{3/4} + n^{1/4} \sd_{2,n}(0) \sP_n(0)^2 \cs(0)^2/T^{1/2}
	\\
	\text{ and }\sa_{5,n}=(\sd_{5,n}(0) + 3 \Paren*{\sd_{3,n}(0)}^2) \sP_n(0)^5 \cs(0)^5/T^{5/4} 
		\\
	+ (5 n^{1/4} \sd_{2,n}(0)  \sd_{3,n}(0) + n^{1/4} \sd_{4,n}(0)) \sP_n(0)^4 \cs(0)^4/T 
		\\
	+ 2 n^{1/2} \Paren*{\sd_{2,n}(0)}^2 \sP_n(0)^3 \cs(0)^3/T^{3/4}. 
 	\end{gather*}
\end{remark}
\begin{remark}
As $\sd_{k,n}$ and $\sP_n$ depend on $\theta$, they cannot be computed
under the true parameter in the context of estimation
(actually, $\sP_n(\theta_n)=0=\sd_{0,n}(\theta_n)$ for the MLE $\theta_n$)
but can be used for statistical hypothesis testing. 
\end{remark}

The result shows a sort of \textquote{phase transition} between $\theta=0$ and $\theta\neq0$.
This is clearly related to the dichotomy in the convergence of $\sd_{2m,n}(0)$ (vanishing) 
and~$n^{1/4} \sd_{2m,n}(0)$ for $m\geq 1$: 
%which is a consequence of the dichotomy exposed in Proposition~\ref{prop:score} :
Observe, for instance, that the second term in equation~\eqref{eq:mle:0} is of order~$1/n^{3/4}$ and in equation~\eqref{eq:mle:neq0} is of order~$1/n^{1/2}$. This latter term,
% of order two, 
in the case~$\theta=0$ goes to 0 as $n$ increases.

Theorem~\ref{thm:asym} and Proposition~\ref{prop:pn_dn}
prove Theorem~\ref{thm:mixed_normality}.
Moreover they imply, for all $m\in \mathbb N$, approximations of the MLE given by the formal power series 
\begin{equation} 
    \label{eq:mle:Mapprox:sd}
	{\theta}_n^{[m]} := \theta+\Phi^{[m]}(\sD_{\cdot,n}(\theta),\sd_{0,n}(\theta)) \quad \text{and} \quad 
	\check{\theta}_n^{[m]}
    := 
	\theta+\Phi^{[m]}(\sD_\cdot(\theta),\sd_{0,n}(\theta)), 
\end{equation}
where $\sD_k(\theta)$ is defined similarly to $\sD_{k,n}(\theta)$ in~\eqref{eq:sD:recursive}
with $\sd_{k,n}$ replaced by $\sd_{n}$ in \eqref{eq:def:dn}.
More precisely,
$\sD_{1}(\theta):=1$ and for $q =2,\ldots, N$, $\sD_q$ is given by the recursive formula:
 \begin{equation} 
     \label{eq:sD:limit:recursive}
	\sD_{q}(\theta) :=  \sum_{m=2}^{q} \sd_{m}(\theta) \sum_{k_1+\dotsb+k_m=q}
    \sD_{k_1}(\theta)\dotsb \sD_{k_m}(\theta).
\end{equation}
Then, for instance,
\begin{equation*}
    \check{\theta}_n^{[3]}
	= \theta
    +
    \sd_{0,n}
    -
    \frac{1}{2}\frac{\xi_2(\theta)}{\xi_1(\theta)}
    \Paren*{ \sd_{0,n}}^2
    + \frac{1}{2}
    \Paren*{
	\Paren*{\frac{\xi_2(\theta)}{\xi_1(\theta)}}^2
    -\frac{\xi_3(\theta)}{3\xi_1(\theta)}}
\Paren*{ \sd_{0,n}}^3.
    %+\grandO\Paren*{\frac{1}{n}}
\end{equation*}
%\end{remark}
%
%
The power series $\theta_n^{[m]}$ 
is the $m$-th order truncation of~\eqref{eq:mle:neq0}.
It has random coefficients~$\sD_{\cdot,n}(\theta)$ 
and random argument~$\sd_{0,n}(\theta)$
involving the score and its derivatives.
The proxy~$\check{\theta}_n^{[m]}$ 
is a power series with deterministic coefficients $\sD_{\cdot}(\theta)$ (limit of $\sD_{\cdot,n}(\theta)$) 
and the same random argument $\sd_{0,n}(\theta)$ related to the score and its first derivative.
Both $\theta_n^{[m]}$ and $\check{\theta}_n^{[m]}$ are non-linear expressions of $\sd_{0,n}(\theta)$.

The following result enlightens the coefficients' behavior as $\theta$ varies in $(-1,1)$.
It is proved in Section~\ref{sec:score_coefficients}.

\begin{lemma} 
    \label{lem:xi_summary}
	Let $\Set{\xi_m}_{m\in \NN}$ defined in~\eqref{eq:xi:def}.
	Then for all $m\in \NN$:
	\begin{enumerate}[thm]
		\item \label{item:xi:sum:1}
			$\xi_0 \equiv 0$;
		\item \label{item:xi:sum:2}
			$\theta \mapsto \xi_{2m+1}(\theta)$ is even;
		\item \label{item:xi:sum:3}
			There exists $-\infty <c_{m,1} < c_{m,2}<0$ such that $-c_{m,1} \leq (1-\theta^2)^{2 m+1} \xi_{2m+1}(\theta)\leq c_{m,2}$ for all $\theta \in (-1,1)$;
		\item \label{item:xi:sum:4}
			$\xi_{2m}(\theta)$ is odd, in particular $\xi_{2m}(0)=0$;
		\item \label{item:xi:sum:5}
			There exists $c_{2m}\in (0,\infty)$ such that $-c_{2m} \leq (1-\theta^2)^{2 m} \xi_{2m}(\theta)<0$ for all~$\theta \in (0,1)$.
	\end{enumerate}
\end{lemma}

The following result is proved by induction on the recursive formula~\eqref{eq:sD:recursive}.

\begin{corollary}[of Proposition~\ref{prop:pn_dn} and Lemma~\ref{lem:xi_summary}]
    \label{cor:sD}
It holds that 
\begin{equation*}
    \sD_{k,n}(\theta)
    \cvproba[n\to\infty] \sD_k(\theta).
\end{equation*}
Besides, for any $k\geq 1$, 
\begin{equation*}
    \abs{\sD_k(\theta)\cs(\theta)^k}\leq \frac{C_k}{(1-\theta^2)^{k/2-1}}
\end{equation*}
for a constant $C_k\geq 0$.
\end{corollary}

\begin{figure}
    \begin{center}
	\includegraphics{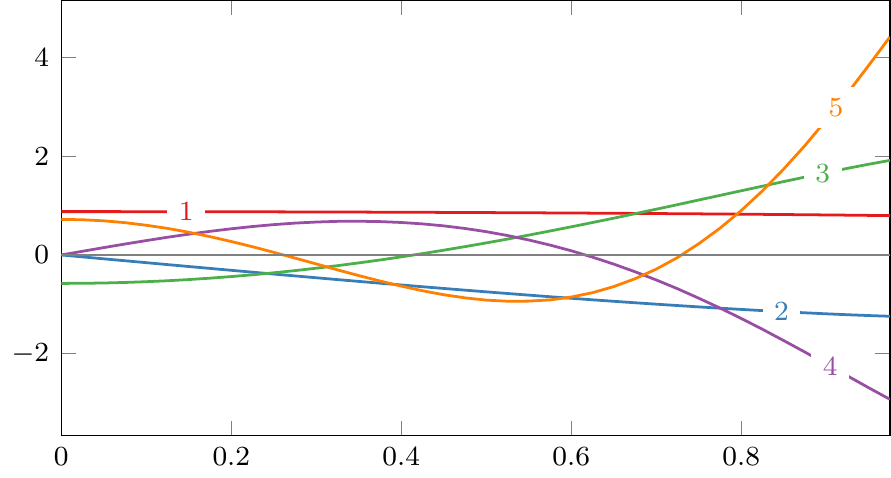}
	\caption{\label{fig:sD:coeff} Functions $\theta\mapsto \sD_{k}(\theta)\cs(\theta)^k(1-\theta^2)^{k/2-1}$ 
	    with $\sD_{k}$ defined by~\eqref{eq:sD:limit:recursive}
    for $\theta$ restricted to $[0,1]$ and $m=1,\dotsc,5$.}
\end{center}
\end{figure}

We plot the coefficients $\theta\mapsto \sD_k(\theta)\cs(\theta)^k$ in Figure~\ref{fig:sD:coeff}, 
up to this scaling factor $\theta\mapsto (1-\theta^2)^{1-k/2}$.

\begin{remark} \label{rem:sD}
	The coefficients $\sD_{k,n}(\theta)$ converge towards their limits at rate of $n^{1/4}$. 
	%This follows from the expression of the coefficients~\eqref{eq:sD:recursive} in terms of the coefficients $\sd_{k,n}(\theta)$
	%which converge towards their limits at rate of $n^{1/4}$ (see Proposition~\ref{prop:pn_dn}).
	This is a consequence of definition~\eqref{eq:sD:recursive} and Proposition~\ref{prop:pn_dn}.
\end{remark}

\begin{figure}
    \begin{center}
    \includegraphics{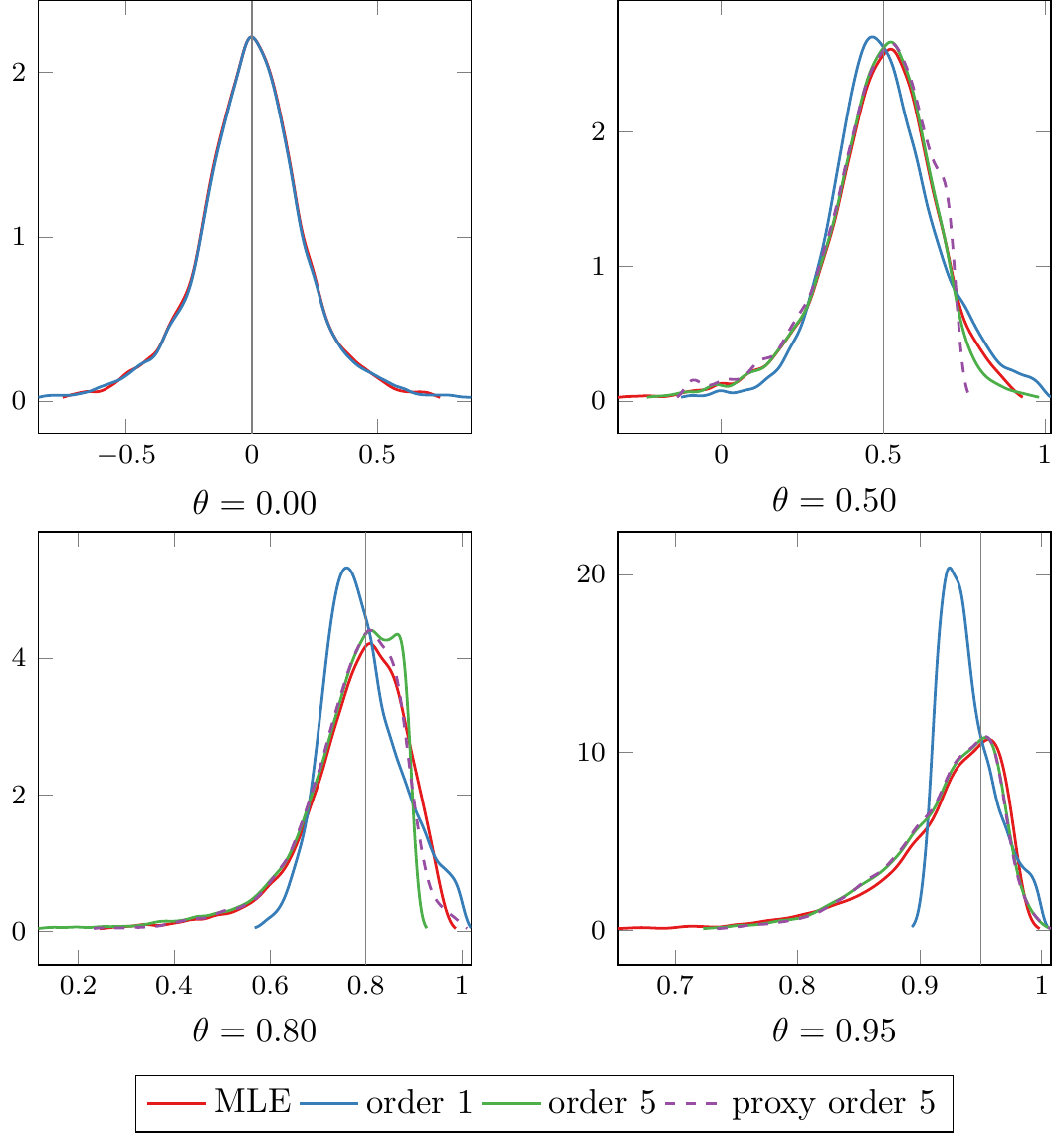}
    \caption{\label{fig:approximation:mle} {Density of MLE, $\theta_n^{[1]}=\theta+\sd_{0,n}(\theta)$, and the fifth order expansions~$\theta_n^{[5]}$ and~$\check{\theta}_n^{[5]}$~\eqref{eq:mle:Mapprox:sd}, 
for several values of $\theta$ using $10\,000$
samples of the SBM using $\Delta t=10^{-3}$ for $T=1$ (thus, $n=1000$).}}
\end{center}
\end{figure}

In Figure~\ref{fig:approximation:mle} we show the empirical densities of the MLE $\theta_n$, 
$\theta_n^{[1]}
%=\theta+\Phi^{[1]}(\sD_{\cdot,n}(\theta),\sd_{0,n}(\theta))
= \theta + \sd_{0,n}(\theta)$,  
$\theta_n^{[5]}$, %=\theta+ \Phi^{[5]}(\sD_{\cdot,n}(\theta),\sd_{0,n}(\theta))$, 
and of the proxy %truncated at fifth order 
$\check{\theta}_n^{[5]}$%=\Phi^{[5]}(\sD_{\cdot}(\theta),\sd_{0,n}(\theta))$.
(see~\eqref{eq:mle:Mapprox:sd} for the definitions of $\theta_n^{[5]}$ and $\check{\theta}_n^{[5]}$).
For $\theta=0$, $\sd_{0,n}(0)$ replicates already very well the MLE behavior, therefore there $\theta_n^{[5]}$ and $\check{\theta}_n^{[5]}$ are not plotted.
For $\theta \neq 0$, we observe that, 
while the lower order expansions replicate worse the MLE behavior,
$\theta_n^{[5]}$ and $\check{\theta}_n^{[5]}$ do it quite well, and so do the higher order expansions.
The expansions $\theta_n^{[5]}$ and $\check{\theta}_n^{[5]}$ are close one to another. Indeed they are respectively random and deterministic polynomials of $\sd_{0,n}(\theta)$ such that, by Corollary~\ref{cor:sD}, the random coefficients converge to the deterministic ones.
To introduce the next paragraph, observe that the density of $\theta_n^{[1]}=\theta + \sd_{0,n}(\theta)$ is skewed to the left when~$\theta$ is close to $1$.
Of course, since their density is close to the one of the MLE's, $\theta_n^{[5]}$ and $\check{\theta}_n^{[5]}$ exhibit a skewed empirical distribution.
\begin{remark}
Let us remind that as for $\theta_n-\theta$, $n^{1/4}\sd_{0,n}(\theta)$ has asymptotic symmetric distribution: Proposition~\ref{prop:pn_dn} ensures that
$n^{1/4}\sd_{0,n}(\theta)=\sP_n(\theta) \cs(\theta) T^{-1/4}$
is asymptotically $\cL$-mixed normal distributed.
Therefore the alternative proxys $\Phi^{[m]}(\sD_{\cdot}(\theta),\sP_\infty(\theta) \cs(\theta) (n T)^{-1/4})$
where $\sP_\infty$ follows a $\cL$-mixed normal distribution, would not replicate skewness and they would need higher order expansions than the proxys $\check{\theta}^{[m]}_n$ in order to have empirical distribution which gets close to the MLE's one.
\end{remark}
%%

%%%%
\paragraph*{Skewness.}
We have observed skewness for the empirical distribution $\sd_{0,n}(\theta)$ and of the MLE.
This is due to the fact that the empirical density of the score shows that the more $\theta$ is close to $1$ the more the score is skewed to 
the left. 
%This induces the same behavior on the MLE. 
However, the correlation between MLE and score is non-trivial : non-linear dependence. 
The approximations of the MLE in~\eqref{eq:mle:Mapprox:sd} 
are non-linear functions of the score and its derivative, actually of $\sd_{0,n}$ which is basically their ratio.
Since $\sD_2(\theta) =\sd_2(\theta)=-\frac12\frac{\xi_2(\theta)}{\xi_1(\theta)}$, 
Lemma~\ref{lem:xi_summary} establishes that $\sgn(\sD_2(\theta))=-\sgn(\theta)$.
This favors a skewness to the left (resp. right)
when $\theta>0$ (resp. $\theta<0$) of the distribution of the proxys~\eqref{eq:mle:Mapprox:sd}.
In other words, the approximations of the MLE proposed above capture the skewness of the~MLE.
%Although the expansions suggest the skewness, of course we could conclude from this that the MLE is skewed though.

\paragraph*{Boundary layer effect.}
Let us set 
\begin{equation}
    \label{eq:f:def}
    f_m(\theta):=\frac{(-1)^m (1-\theta^2)^m}{m!} \xi_m(\theta).
\end{equation}

\begin{figure}
    \begin{center}
    \includegraphics{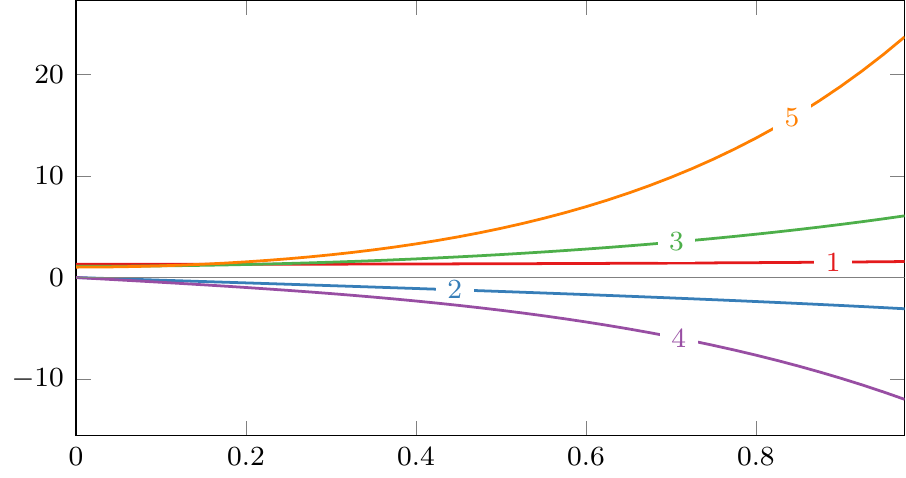}
    \caption{\label{fig:f:coeff} Functions $\theta\mapsto f_m(\theta)$ 
	defined by \eqref{eq:f:def} for 
    $\theta$ restricted to $[0,1]$ and $m=1,\dotsc,5$.}
\end{center}
\end{figure}

Then $f_{2k+1}$ is positive and bounded (from above and below by a positive constant), $f_{2k}(0)=0$ and $f_{2k}(\theta)$ is bounded.
Then, \eqref{eq:mle:Mapprox:sd} rewrites:
\begin{multline*}
	\check{\theta}^{[3]}_n 
	     = \theta 
	    + \frac{ \sqrt{1-\theta^2} }{\sqrt{f_1(\theta)}} 
	    \frac{\sP_n(\theta)}{(Tn)^{1/4}}
	    +
	    \frac{f_2(\theta)}{f_1(\theta)^2} \Paren*{\frac{\sP_n(\theta)}{(Tn)^{1/4}}}^2
    \\
	    + \frac{1}{\sqrt{1-\theta^2}} \frac{ 2 (f_2(\theta))^2 - f_3(\theta) f_1(\theta)}{f_1(\theta)^{7/2}}\Paren*{\frac{\sP_n(\theta)}{(Tn)^{1/4}}}^3,
    \end{multline*}
    that is 
\begin{align*}
    \sD_1(\theta)\cs(\theta)&= \frac{ \sqrt{1-\theta^2} }{\sqrt{f_1(\theta)}},
    \quad
    \sD_2(\theta)\cs(\theta)^2=\frac{f_2(\theta)}{f_1(\theta)^2}
    \\
    \text{ and }
\sD_3(\theta)\cs(\theta)^3&=\frac{1}{\sqrt{1-\theta^2}} \frac{ 2 f_2(\theta)^2 - f_3(\theta) f_1(\theta)}{f_1(\theta)^{7/2}}.
\end{align*}
The term of order $2$ vanishes for $\theta=0$ as $f_2(0)=0$. 
With this expansion, one sees a \textquote{boundary layer} effect for $\abs{\theta}$ close to $1$ in the explosion of the third order coefficient.
Indeed if we push the expansion up to order $k$ (here we stopped at~$k=3$), the corresponding coefficient explodes as $(1-\theta^2)^{(k-2)/2}$ when $\abs{\theta}$ is close to~$1$.
We have performed numerical simulations that suggest that no cancellation effect occurs
so that the approximation by a polynomial expansion is no longer suitable.

In~\cite{lm22} we discussed a similar boundary layer phenomenon in the case of the Binomial family.
	One could expect the same to happen for the skewness parameter~$\theta$ 
	because of the pathwise construction of SBM done 
	associating independent Bernoulli random variables with parameter~$(1+\theta)/2$
	to the excursions from 0 of a reflected Brownian motion 
	and flipping each excursion based on the result of the Bernoulli random variable.

%%%%
\paragraph*{Change of variable and change of coordinates.} Combining the Faà di Bruno 
formula with \eqref{eq:mle:neq0}, 
one may given some explicit expansion of $\varphi(\theta_n)$ for 
any analytic function $\varphi$. Similarly, one may also consider $\theta_n$
in another system of coordinates, as discussed in \cite{lm22}.
On that point, two changes of variables appear to be natural:~$\varphi_1(\theta)=\sqrt{1-\theta^2}$ and $\varphi_1(\theta)=1/\sqrt{1-\theta^2}$. The
latter one stabilizes the variance.
However, no change of coordinate impacts the asymptotic behavior of the score.
Besides, we found through numerical experiments that using~$\varphi_1$ or~$\varphi_2$ does not improve
Wald confidence intervals.

%%%%%%%%%%%%%%%%%%%%%%%%%%%%%%%%%%%%%%%%%%%%%%%%%%%%%%%%%%%%%%%%%%%%%%%
%\subsection{Outline of the rest of the paper}
%We have stated several results that we prove in the next sections.
%%
%In Section~\ref{sec:score}, we study the asymptotic behavior of the score and its derivative and we prove Proposition~\ref{prop:pn_dn}.
%The results, Propositions~\ref{thm:rate:conv} and~\ref{prop:score}, are interesting on their own
%%, and they are essential to our proofs of the asymptotic results proposed above, in particular of the key Proposition~\ref{prop:pn_dn}. The proofs of Propositions~\ref{prop:pn_dn},\ref{thm:rate:conv}, and~\ref{prop:score} 
%and their proofs are given in Section~\ref{sec:proofs}.
%%
%In Section~\ref{sec:score_coefficients}, we prove Proposition~\ref{prop:s_theta} and Lemma~\ref{lem:xi_summary}.
%We actually establish more precise statements together additional properties useful in numerical studies and applications. 
%For instance, in Section~\ref{sec:coeff_expansion}, we provide an expansion of $\theta \mapsto \xi_m(\theta)$ around 0.
%%
%To conclude, in Section~\ref{sec:numerical}, numerical studies and conclusions are proposed.

%%%%%%%%%%%%%%%%%%%%%%%%%%%%%%%%%%%%%%%%%%%%%%%%%%%%%%%%%%%%%%%%%%%%%%
%%%%%%%%%%%%%%%%%%%%%%%%%%%%%%%%%%%%%%%%%%%%%%%%%%%%%%%%%%%%%%%%%%%%%%
\section{Asymptotic behavior of the score and its derivatives}
\label{sec:score}

In this section, we provide results which are necessary to 
prove Proposition~\ref{prop:pn_dn}.
These results are an application of the results in 
\cite{jacod98} for the Brownian motion,
and on the ones from~\cite{lejay2019} (convergence) and~\cite{mazzonetto19a} (Central Limit Theorem) 
for the (true) SBM. 
%A related result may also be found in~\cite{robert_2022}.

Let us recall that $\coef{S}_{m}$ in~\eqref{eq:coefS} is nothing else that a rescaled derivative of the score
and $\xi_m(\theta)$ is a constant defined in~\eqref{eq:xi:def}.
By Lemma~\ref{lem:xi_summary}, it holds that for all $\theta \in (-1,1)$,
$\xi_0(\theta)=0$ and $\xi_1(\theta)$ is negative, and that for all integer $m$, $\xi_{2m}(0)=0$.

\begin{proposition}
    \label{prop:score}
    For any $\theta\in(-1,1)$, 
	%there exists a family $\Set{\xi_k(\theta)}_{k\geq 0}$ defined by~\eqref{eq:xi} below 
	the family $\Set{\xi_k(\theta)}_{k\geq 0}$ defined by~\eqref{eq:xi:def} is
    such that: 
    \begin{enumerate}[thm]
\item \label{prop:score:ii}
     Under $\PP_\theta$, 
    there exists a standard (with mean $0$ and variance $1$) Gaussian random variable $H$ independent from $\cF_T$
    (the probability space $(\Omega,\cF,\PP_0)$ has been extended to carry $H$)
    such that for any $m\geq 0$, 
    \begin{align*}
	& n^{1/4}\coef{S}_0(n,\theta) \cvstable[n\to\infty] \sqrt{-\xi_1(\theta)}\sqrt{L_T}H,
	\\
	&(\coef{S}_1(n,\theta),\dotsc,\coef{S}_m(n,\theta))	
	\cvproba[n\to\infty] 
	    (\xi_1(\theta)L_T,\dotsc,\xi_m(\theta)L_T).
    \end{align*}
    Thanks to the property of the stable convergence, 
    we have joint stable convergence of $(n^{1/4}\coef{S}_0(n,\theta),\coef{S}_1(n,\theta),\dotsc,
\coef{S}_m(n,\theta))$.
\item \label{prop:score:i}
    For $\theta=0$ (the SBM is actually a Brownian motion), under $\PP_0$, 
    there exists a Gaussian family $H=\Set{H_k}_{k=0,\dotsc,m}$ independent from $\cF_T$
    with mean $0$ and covariance $\Cov(H_j,H_\ell)=\Psi_{2j,2\ell}$ for $0\leq j,\ell\leq m$ 
    described in Section~\ref{sec:Psi} 
    (the probability space $(\Omega,\cF,\PP_0)$ has been extended to carry $H$)
    such that for any $m\geq 0$,
	\begin{align*}
	& \Set{ n^{1/4} \coef{S}_{2k}(n,0)}_{k=0,\dotsc,m} 
	    	\cvstable[n\to\infty] 
		\Set{\sqrt{L_T}H_{k}}_{k=0,\dotsc,m},
	\\
	&
	\Set{\coef{S}_{2k+1}(n,0)}_{k=0,\dotsc,m} 
		\cvproba[n\to\infty] 
		\Set{\xi_{2k+1}(0) L_T}_{k=0,\dotsc,m}.
    \end{align*}
    Thanks to the property of the stable convergence, 
    we have joint stable convergence of 
    $(n^{1/4}\coef{S}_0(n,\theta),\coef{S}_1(n,\theta),
    n^{1/4}\coef{S}_{2}(n,\theta),\dotsc)$.
    \end{enumerate} 
\end{proposition}

\begin{remark}
    \label{rem:score}
    Following \cite{jacod98,mazzonetto19a},
a stronger result holds:  on an extended probability space,
there exist a $m$-dimensional Brownian motion  
$\Set{W_k}_{k\geq 0}$ independent from~$L$ and a $m\times m$-matrix $\Xi(\theta)$ 
such that 
\begin{multline}
    \label{eq:rate:conv}
    \Set*{\frac{1}{n^{1/4}}
	\Paren*{
	    \frac{1}{\sqrt{n}}
	    \sum_{i=1}^{\Floor{nt}}
	    \partial^\ell k_\theta(X_{i-1}\sqrt{n},X_i\sqrt{n})
	-\xi_{\ell}(\theta) L_t}
    }_{\ell=0,\dotsc,m,\ t\in[0,T]}
    \\
    \cvstable[n\to\infty] \Set{\Xi(\theta) W(L_t)}_{t\in [0,T]}.
\end{multline}
A closed-form --- yet cumbersome --- expression exists for the matrix $\Xi(\theta)=\left( \Xi_{i,j}(\theta)\right)_{i,j=0,\ldots,m}$.
Here, we focus on the particular cases $m=0$ and $\theta=0$ 
%which are of interest for statistical applications 
and we study $\Set{\xi_k(\theta)}_{k=1,\dotsc,\max(m,1)}$ and $\Xi$
for these cases.
\end{remark}
\begin{remark} \label{rem:rate:conv}
	When $\xi_m(\theta)=0$, \eqref{eq:rate:conv} rewrites
	\begin{equation*}
		n^{1/4}	\coef{S}_m(n,\theta)
		\cvstable[n\to\infty]  \Xi_{m,m}(\theta) \sqrt{L_T}G
	\end{equation*}
	for some constant $\Xi_{m,m}(\theta)$, and $G\sim\cN(0,1)$ independent of $L_T$.
	In particular,
	\begin{enumerate}[thm]
	    \item For any $\theta\in(-1,1)$ and $m=0$, 
	    \begin{equation*} 
		\Xi_{0,0}(\theta)
		=\sqrt{-\xi_1(\theta)}
		=\frac{1}{\cs(\theta)}.
	    \end{equation*}
	    This follows from \cite{mazzonetto19a} and Lemma~\ref{lem:chi:m=1}. 
%		Indeed there a more complicate expansion is given involving the term
%	    \begin{equation}
%		\label{eq:average:0}
%		\EE_\theta\Prb{k_\theta(X_0,X_1)\given X_0=x}
%		=
%		\int_{\RR} \frac{\partial_\partial p_\theta(1,x,y)}{p_\theta(1,x,y)}
%		p_\theta(1,x,y)\vd y=0,\ \forall x\in\RR
%	    \end{equation}
%	    (see Lemma~\ref{lem:chi:m=1} for the latter equality showing $\xi_0(\theta)=0$). 
%		Since the latter term is 0, the expression given in~\cite{mazzonetto19a} simplifies to the one given above.
    \item For $\theta=0$, 
	$\Xi(0)\Xi(0)^\mathrm{T}=\Psi$ for $\Psi$ appearing
	in Proposition~\ref{prop:score}.\ref{prop:score:i} and discussed in Section~\ref{sec:Psi}.
\end{enumerate}
\end{remark}
%%

%%%%%%%%%%%%%%%%%%%%%%%%%%%%%%%%%%%%%%%%%%%%%%%%%%%%%%%%%%%%%%%%%%%%%%

\paragraph*{Statistical implications.}
The convergence results of Proposition~\ref{prop:score} are the key for
the convergence of the MLE in Theorem~\ref{thm:mixed_normality}. 
Besides, we could also use these results to construct
\begin{itemize}[noitemsep,topsep=-\parskip]
    \item estimators of the local time using $\coef{S}_1(n,\theta_n)/\xi_1(\theta_n)$;
    \item Wald confidence interval using $\cs(\theta_n)$ as a substitute for $\cs(\theta)$;
    \item hypothesis testing on the true value of $\theta$ using 
	either a confidence interval around $\theta$, or $\coef{S}_0(n,\theta)^2/\coef{S}_1(n,\theta)$
	which behaves asymptotically as a $\chi^2_1$ distribution. 
\end{itemize}

%%%%%%%%%%%%%%%%%%%%%%%%%%%%%%%%%%%%%%%%%%%%%%%%%%%%%%%%%%%%%%%%%%%%%%
\subsection{Rate of convergence}

\label{sec:rate}

In this section, 
we estimate numerically the rate of convergence towards the the limit distribution
of the score.
More precisely, we study empirically the rate convergence 
of $\coef{S}_k(n,\theta)$ in~\eqref{eq:coefS} towards its limit. 
In particular, 
we show that the speeds deteriorates as $\abs{\theta}$ becomes
close to $1$.

\begin{figure}
    \begin{center}
    \includegraphics{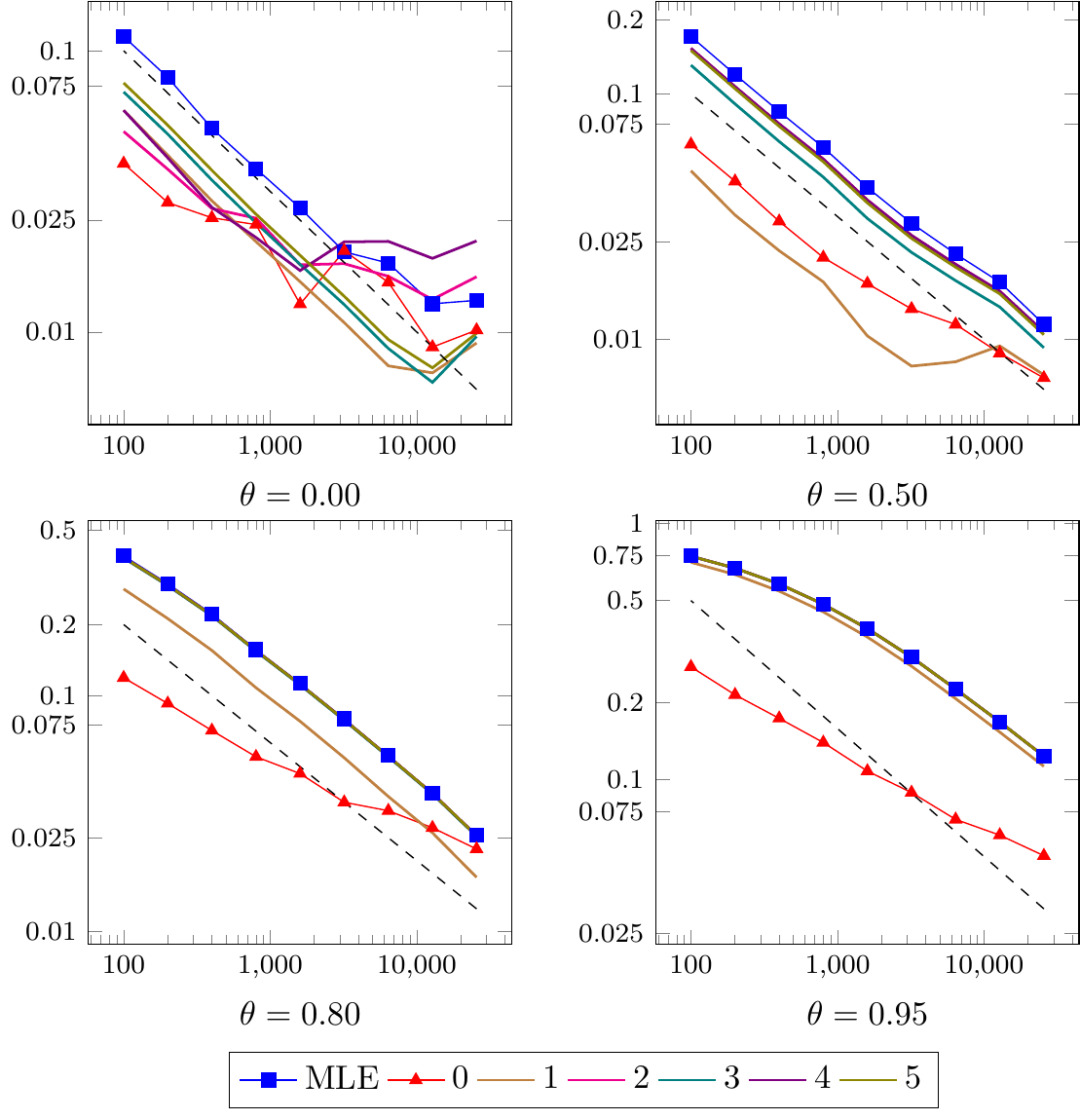}
    \caption{\label{fig:rate} Behavior of $n\mapsto \Delta_{\mathrm{KS}}(n)$ using $10\,000$
    samples of the SBM with $T=1$, and comparison with $C/\sqrt{n}$ (dashed line), in log-log scale, 
    for statistics \ref{stat:A} (MLE), \ref{stat:B} ($0$ for any $\theta$ 
    and $2,4$ if $\theta=0$) and \ref{stat:C} in the other cases.}
\end{center}
\end{figure}

We plot in Figure~\ref{fig:rate}
the Kolmogorov-Smirnov distance~$\Delta_{\mathrm{KS}}(n)$
between the empirical distributions 
\begin{enumerate}[noitemsep,topsep=-\parskip,label={(\alph*)}]
    \item
	\label{stat:A} 
	of $\cs(\theta)^{-1}n^{1/4}(\theta_n-\theta)$ and $H/\sqrt{L_T}$
	with $H\sim\cN(0,1)$, independent from the local time.
    \item\label{stat:B} of $\coef{S}_k(n,\theta)/\Xi_{k,k}(\theta)$ and $\sqrt{L_T}H$ with $H\sim\cN(0,1)$, independent
	from the local time, for $k=0$ or for $\theta=0$ and $k$ even.
    \item\label{stat:C} of $\coef{S}_k(n,\theta)/\xi_k(\theta)$ and $L_T$ in other cases.
\end{enumerate}

Even for small moderate values of 
the sample size $n$, the asymptotic regime is reached
most of the case. The 5\,\%-quantile of the Kolmogorov-Smirnov
statistics is~$0.52$ while the~1\,\%-quantile
is~$0.44$. 

In each case, the statistics $\Delta_{\mathrm{KS}}(n)$
behaves like $C/n^\eta$, with $\eta$ close to $1/2$. 

For the score of order $0$ or when $\theta$ is close to $1$, the rate $\eta$ is
lower than $1/2$.  The constant $C$ is greater for the MLE than for the other
Kolmogorov-Smirnov statistics.  In each case, the constant increases with
$\theta$, as well as the ratio between the constant $C$ for the MLE and the one
for the score of order $0$.

Note that we are close to the Berry-Esseen rate 
even though our sample is not composed of independent identically distributed random variables. 

In Figure~\ref{fig:rate}, when the Kolmogorov-Smirnov distances are small
and~$n$ is big, the plotted distance is affected by the Monte-Carlo error. 
For $\theta=0$, the terms of order~$2$ and~$4$ in the expansion of the
MLE~\eqref{eq:mle:neq0} vanish.  The asymptotic law of these coefficients is
the one given in~\ref{stat:B}.  Note that, in the alternative expansion of the
MLE~\eqref{eq:mle:0}, the terms of odd order of the latter expansion involve
the statistics in~\ref{stat:B} and~\ref{stat:C}.  For $\theta\neq0$, the effect
of the terms of orders $k\geq 2$ in the expansion~\eqref{eq:mle:neq0}, which
behave asymptotically as $(1-\theta^2)^{1-k/2}$, have some increasing effect on
the~MLE.

%%%%%%%%%%%%%%%%%%%%%%%%%%%%%%%%%%%%%%%%%%%%%%%%%%%%%%%%%%%%%%%%%%%%%%
%%%%%%%%%%%%%%%%%%%%%%%%%%%%%%%%%%%%%%%%%%%%%%%%%%%%%%%%%%%%%%%%%%%%%%
\section{The limiting coefficients}
\label{sec:score_coefficients}

Our aim is to study the limiting coefficients $\cs(\theta)$
in~\eqref{eq:s_theta:def}, $\xi_m(\theta)$ in~\eqref{eq:xi:def} and $\Psi$
appearing in Proposition~\ref{prop:score}.  We also prove
Proposition~\ref{prop:s_theta} and Lemma~\ref{lem:xi_summary}.  We actually
establish additional properties useful in numerical studies and applications.
For instance, in Section~\ref{sec:coeff_expansion}, we provide an expansion of
$\theta \mapsto \xi_m(\theta)$ around 0.  In Section~\ref{sec:Psi},~we study
the matrix $\Psi$ in Proposition~\ref{prop:score}.\ref{prop:score:i}.

%%%%%%%%%%%%%%%%%%%%%%%%%%%%%%%%%%%%%%%%%%%%%%%%%%%%%%%%%%%%%%%%%%%%%%

\subsection{The coefficients \texorpdfstring{$\xi_k$}{ξₖ}}

We study up to Section~\ref{sec:coeff_expansion} the limiting coefficients $\xi_k$.

We take profit from the explicit expression of the density, hence of $k_\theta$.
We conveniently rewrite 
\begin{equation*}
    k_\theta(x,y):=
    \begin{cases}
	\dfrac{1}{\theta+\exp(2xy)}&\text{ if }x\geq 0,\ y\geq 0,\\
	\dfrac{1}{\exp(2xy)-\theta}&\text{ if }x\leq 0,\ y\leq 0,\\
	\dfrac{1}{\theta-1}&\text{ if }x\geq 0,\ y\leq 0,\\
	\dfrac{1}{1+\theta}&\text{ if }x\leq 0,\ y\geq 0.
    \end{cases}
\end{equation*}

%%%
\begin{lemma}
    \label{lem:k}
For any $m\geq 0$, 
\begin{equation}
    \label{eq:derk}
    \partial_\theta^m k_\theta(x,y)=m!(-1)^m k_\theta^{m+1}(x,y) 
\end{equation}
for any $x,y\in\RR$.
\end{lemma}
%%%

%%
\begin{proof}
    Let us note first that 
\begin{multline}
    \label{eq:derk2}
    \partial_\theta k_\theta(x,y)=\frac{-\sgn(y)^2}{(\sgn(y)\theta+\exp(2(xy)^+))^2}
    =\frac{-1}{(\sgn(y)\theta+\exp(2(xy)^+))^2}
    \\
    =-k_\theta(x,y)^2.
\end{multline}
Thus, \eqref{eq:derk} is true for $m=0,1$. Let us assume that \eqref{eq:derk} is
true for some $m\geq 1$. With \eqref{eq:derk2}, 
\begin{multline*}
    \partial_\theta^{m+1}k_\theta(x,y)
    = (-1)^m m! \partial_\theta(k_\theta(x,y)^{m+1})
    \\
    =(-1)^m (m+1)!\partial_\theta k_\theta(x,y)\cdot k_\theta(x,y)^m
    =(-1)^{m+1} (m+1)! k_\theta(x,y)^{m+2}.
\end{multline*}
This proves that \eqref{eq:derk} is true for $m+1$, and thus for any $m$.
\end{proof}

\begin{notation}
    \label{not:chi}
We define for some integer $m\geq 1$, 
\begin{align}
    \label{eq:chi:1}
    \chi^{\pm\pm}_m(\theta):=\iint_{\RR^2} \ind{\pm x\geq 0,\pm y\geq 0} k_\theta(x,y)^m p_\theta(1,x,y)\mu_{\theta}(x)\vd x\vd y,
    \\
		\label{eq:chi:decomp}
    \text{and }
    \chi_m(\theta):=\chi_m^{++}(\theta) + \chi_m^{--}(\theta)+\chi^{+-}_m(\theta)+\chi^{-+}_m(\theta).
\end{align}
\end{notation}

With Lemma~\ref{lem:k}, the coefficients $\xi_m$'s given by \eqref{eq:xi:def} are 
related to the $\chi_m$'s by 
\begin{equation}
    \label{eq:xi}
    \xi_m(\theta)=m!(-1)^m\chi_{m+1}(\theta).
\end{equation}

We now study $\chi_m(\theta)$ for $\theta\in(-1,1)$ and various values of $m$.
 
We note the relations
\begin{align}
    \label{eq:sym:1}
    p_\theta(1,-x,-y)&=p_{-\theta}(1,x,y),\\
    \label{eq:sym:2}
    k_\theta(-x,-y)&=-k_{-\theta}(x,y)\\
    \label{eq:sym:3}
    \text{and }\mu_\theta(-x)&=\mu_{-\theta}(x). 
\end{align}

A direct consequence of the symmetry relations \eqref{eq:sym:1}-\eqref{eq:sym:3} is the following
lemma

%%%
\begin{lemma}[Symmetry relation]
    \label{lem:chi:symm}
    For any $m\geq 1$ and any $\theta\in(-1,1)$, 
\begin{gather}
    \label{eq:sym:4}
    \chi^{--}_m(\theta)=(-1)^m\chi^{++}_m(-\theta)
    \text{ and }
    \chi^{-+}_m(\theta)=(-1)^m\chi^{+-}_m(-\theta).
\end{gather}
\end{lemma}
%%%

%%
\begin{corollary}
    \label{cor:chi}
    For any $\theta\in(-1,1)$ and any $m\geq 1$,
    \begin{enumerate}[thm]
	\item $\theta\mapsto \chi_{2m}(\theta)$ is even and $\chi_{2m}(\theta)>0$. 
	\item $\theta\mapsto \chi_{2m+1}(\theta)$ is odd; in particular $\chi_{2m+1}(0)=0$
	    (we will see below in Lemma~\ref{lem:k:non_vanishing} that $\chi_{2m+1}(\theta)<0$
	    for any $\theta\in(0,1)$).
    \end{enumerate}
\end{corollary}
%%%

%%%
\begin{proof}
Consequence of Lemma~\ref{lem:chi:symm} and for proving $\chi_{2m}(\theta)>0$ one just observes the 
positivity of the integrands of equation~\eqref{eq:chi:1}.
\end{proof}
%%%

%%%
\begin{lemma}[Values for $m=1$]
    \label{lem:chi:m=1}
    For any $\theta\in(-1,1)$, $\chi_1(\theta)=0$.
\end{lemma}
%%%

%%%
\begin{proof}
    This stems from the very definition of $k_\theta$ and the fact that 
    $\int_{-\infty}^{+\infty} p_\theta(1,x,y)\vd y=1$
    so that $\partial_\theta \int_{-\infty}^{+\infty} p_\theta(1,x,y)\vd y=0$.
	In fact, 
	\[\chi_1(\theta)=\int_{\mathbb{R}} \EE_\theta\Prb{k_\theta(X_0,X_1)\given X_0=x} \mu_{\theta}(x) \vd x\]
	and for all $x\in\RR$
	\begin{equation}
		\label{eq:average:0}
		\EE_\theta\Prb{k_\theta(X_0,X_1)\given X_0=x}
		=
		\int_{\RR} \frac{\partial_\partial p_\theta(1,x,y)}{p_\theta(1,x,y)}
		p_\theta(1,x,y)\vd y=0.
	    \end{equation}
Hence $\chi_1(\theta)=0$.
\end{proof}
%%%

%%%%%%%%%%%%%%%%%%%%%%%%%%%%%%%%%%%%%%%%%%%%%%%%%%%%%%%%%%%%%%%%%%%%%%
%%%%%%%%%%%%%%%%%%%%%%%%%%%%%%%%%%%%%%%%%%%%%%%%%%%%%%%%%%%%%%%%%%%%%%
%%%%%%%%%%%%%%%%%%%%%%%%%%%%%%%%%%%%%%%%%%%%%%%%%%%%%%%%%%%%%%%%%%%%%%

\paragraph*{Computation for $m\geq 2$.}
The case $m\geq 2$ and $\theta>0$ is a bit more cumbersome. 

%%%
\begin{lemma}
    \label{lem:chimp}
For $m\geq 2$, $\theta \in (-1,1)$,
\begin{equation}
\label{eq:chi:6}
    \chi^{+-}_{m}(\theta)
    =\frac{(-1)^m}{\sqrt{2\pi}}\frac{1+\theta}{(1-\theta)^{m-1}}
    \text{ and }
    \chi^{-+}_{m}(\theta)
    =\frac{1}{\sqrt{2\pi}}\frac{1-\theta}{(1+\theta)^{m-1}}.
\end{equation}
\end{lemma}
%%%

%%%
\begin{proof}
First,
\begin{equation*}
    \frac{1}{\sqrt{2\pi}}
    \iint_{\RR^2}\ind{x\geq 0,y\geq 0}\exp\left(\frac{-(x+y)^2}{2}\right)\vd x \vd y
    =
    \int_0^{+\infty} \overline{\Phi}(z)\vd z
    =
    \frac{1}{\sqrt{2\pi}}
\end{equation*}
with $\overline{\Phi}(x):=\frac{1}{\sqrt{2\pi}}\int_x^{+\infty} e^{-y^2/2}\vd y$, 
the complementary distribution function of the unit, centered Gaussian random variable.

 Hence, 
\begin{equation*}
    \chi^{+-}_m(\theta)=\iint_{\RR^2} \ind{x\geq 0,y\leq 0}
    k_\theta(x,y)^mp_\theta(1,x,y)\mu_{\theta}(x)\vd x\vd y
    =\frac{-(1+\theta)}{(1-\theta)^{m-1}}
    \int_0^{+\infty}\overline{\Phi}(x)\vd x.
\end{equation*}
The value of $\chi_m^{-+}$ follows from \eqref{eq:sym:4}.
\end{proof}
%%%

We express 
\begin{multline}
    \label{eq:chi:3}
    \chi^{++}_m(\theta)=
    \iint_{\RR^2} \ind{x\geq 0,y\geq 0}
    k_\theta(x,y)^mp_\theta(1,x,y)\mu_{\theta}(x)\vd x\vd y
    \\
    =
    \iint_{\RR^2} \ind{x\geq 0,y\geq 0}
    \frac{1}{\sqrt{2\pi}} (1+\theta)\frac{\exp\Paren*{\frac{-(x-y)^2}{2}}}{(\theta+\exp(2xy))^m}\vd x\vd y
    \\
    +
    \iint_{\RR^2} \ind{x\geq 0,y\geq 0}
    \frac{1}{\sqrt{2\pi}}\theta(1+\theta)\frac{\exp\Paren*{\frac{-(x+y)^2}{2}}}{(\theta+\exp(2xy))^m}\vd x\vd y\\
    =
    \iint_{\RR^2} \ind{x\geq 0,y\geq 0}
    \frac{1+\theta}{\sqrt{2\pi}} \exp\Paren*{\frac{-(x+y)^2}{2}}\frac{1}{(\theta+\exp(2xy))^{m-1}}\vd x\vd y.
\end{multline}
We could also use the following form 
\begin{equation}
    \label{eq:chi:3bis}
    \chi^{++}_m(\theta)
    =
    \iint_{\RR^2} \ind{x\geq 0,y\geq 0}
    \frac{1+\theta}{\sqrt{2\pi}} \exp\Paren*{\frac{-x^2-y^2}{2}}\frac{\exp(xy)+\theta\exp(-xy)}{(\theta+\exp(2xy))^m}\vd x\vd y
\end{equation}
which is suitable for Monte Carlo simulation as 
\begin{equation} 
    \label{eq:MC}
    \chi^{++}_m(\theta)
    =\frac{\sqrt{2\pi}}{4}(1+\theta)
    \EE\Prb*{\frac{\exp(\abs{G\cdot G'})+\theta\exp(-\abs{G\cdot G'})}{(\theta+\exp(2\abs{G\cdot G'}))^m}}
\end{equation}
for two independent Gaussian, unit, centered random variables $G$ and $G'$.

\begin{remark}
    \label{rem:variance-reduction}
    The considerations in the remainder of this section or in Lemma~\ref{lem:xi_summary} (which we are proving in this section), suggest to consider the equivalent expression 
\begin{equation*} 
    \chi^{++}_m(\theta)
    =\frac{\sqrt{2\pi}}{4(1-\theta)^m}(1+\theta)
    \EE\Prb*{(1-\theta)^m\frac{\exp(\abs{G\cdot G'})+\theta\exp(-\abs{G\cdot G'})}{(\theta+\exp(2\abs{G\cdot G'}))^m}}
\end{equation*}
which reduces the variance of the empirical mean without burdening the computational cost.
\end{remark}

We consider first a close-form formula for the integral involved in $\chi^{++}_{m}(0)$.

\begin{lemma}[Some values at $\theta=0$]
    \label{lem:theta:0}
    Let $m\geq 2$, then
    \begin{equation*}
	\sqrt{2\pi} \chi^{++}_m(0)= \frac{1}{2 \sqrt{m(m-1)}} \log \Paren*{2m-1+ 2\sqrt{m(m-1)}} \in (0,1). 
    \end{equation*}
In particular, when $m$ is large, $\sqrt{2\pi} \chi^{++}_m(0)\sim \log(4m)/2m$.
\end{lemma}

\begin{remark}
    We could also write $\chi^{++}_m(0)=\frac{1}{2\sqrt{m(m-1)}} \atanh\frac{2\sqrt{m (m-1)}}{2m-1}$.
\end{remark}

\begin{proof}[Proof of Lemma~\ref{lem:theta:0}]
From~\eqref{eq:chi:3} we see that $\sqrt{2\pi} \chi^{++}_m(0)= \alpha_{2m-1}$ where
\begin{equation*}
	\alpha_n :=\int_0^{+\infty}\int_0^{+\infty}\exp\Paren*{\frac{-x^2-y^2}{2}- n xy}\vd x\vd y
\end{equation*}
for $n\geq 1$.
Let $n\geq 2$ be fixed.
    Using the change of variable $(x,y)=(\sqrt{v\log(u)},\sqrt{\log(u)/v})$
    for $u\geq 1$, $v\geq 0$, 
    we obtain that 
    \begin{equation*}
	\alpha_n=\int_0^{+\infty}
	\int_1^{+\infty} \frac{1}{2vu}\frac{1}{u^{n+v/2+1/2v}}\vd u \vd v
	=\int_0^{+\infty}
	\frac{1}{v^2+2 n v+1}\vd v.
    \end{equation*}
    % Formula 108 in 1999_Bookmatter_WorkshopCalculusWithGraphingCa.pdf
    Since $n\geq 2$, from (3.3.17) in \cite{abramowitz}, this leads 
\begin{equation*}
   \alpha_n =
	\frac{1}{2 \sqrt{n^2 -1}}
	\log\frac{n + \sqrt{n^2 -1}}{n - \sqrt{n^2 -1}}
	= 	\frac{1}{\sqrt{n^2 -1}}
	\log( n+\sqrt{n^2 -1}) \in (0,1)
\end{equation*}
which can be also written as $\frac1{\sqrt{n^2-1}} \atanh\frac{\sqrt{n^2-1}}{n}$.
\end{proof}

The following inequalities are consequences of~\eqref{eq:chi:3}. For all $\theta\in [0,1)$, $m\geq 2$:
\begin{equation} 
    \label{eq:chi:bounds}
    \pm \frac{1}{\sqrt{2\pi}} \frac{1}{(1\pm\theta)^{m-2}} \sqrt{2\pi}\chi^{++}_{m}(0) 
	\leq
	    \pm \chi^{++}_{m}(\pm\theta) 
	\leq
	    \pm \frac{1\pm\theta}{\sqrt{2\pi}}\sqrt{2\pi}\chi^{++}_{m}(0).
\end{equation}

%%%
\begin{lemma}
    \label{lem:k:non_vanishing}
    For any $m\geq 1$ and $\theta\in (-1,1)$ with $\theta\neq 0$, then $\xi_m(\theta)\neq0$.
\end{lemma}
%%%%

%%%
\begin{proof}
Equality~\eqref{eq:xi} and Corollary~\ref{cor:chi} ensure that it suffices to prove that $ (0,1) \ni \theta \mapsto \chi_{2m+1}(\theta)<0$.
Let $m\geq 1$ and $\theta \in (0,1)$ be fixed.

Equations~\eqref{eq:sym:4}, the explicit expression for $\chi^{-+}_{2m+1}(\theta)$ in~\eqref{eq:chi:6}, the upper inequality in~\eqref{eq:chi:bounds}
and the fact that $\sqrt{2\pi}\chi^{++}_{2m+1}(0) < 1$ (see Lemma~\ref{lem:theta:0}) yield
\begin{equation*}
\begin{split} 
\chi_{2m+1}(\theta) 
	&=\chi^{++}_{2m+1}(\theta) - \chi^{++}_{2m+1}(-\theta) + \chi^{+-}_{2m+1}(\theta) - \chi^{+-}_{2m+1}(-\theta) \\
	& \leq \frac{1}{\sqrt{2\pi}} \left( 2\theta \sqrt{2\pi} \chi^{++}_{2m+1}(0)- \frac{1+\theta}{(1-\theta)^{2m}} + \frac{1-\theta}{(1+\theta)^{2m}} \right) \\
	& <  \frac{1}{\sqrt{2\pi}} \left( 2\theta - \frac{1+\theta}{(1-\theta)^{2m}} + \frac{1-\theta}{(1-\theta)^{2m}} \right)  
	=\frac{2\theta}{\sqrt{2\pi}} \left( 1 - \frac{1}{(1-\theta)^{2m}}\right)<0.
\end{split}
\end{equation*}
This proves that $\xi_m(\theta)$ does not vanish for $\theta\neq0$.
\end{proof}
%%%

We are now ready to provide the proof of Lemma~\ref{lem:xi_summary}.
The proof of Proposition~\ref{prop:s_theta} is also provided and requires studying $\chi_2(\theta)=-\xi_1(\theta)=\cs(\theta)^{-2}$.

%%%%%%%%%%%%%%%%%%%%%%%%%%%%%%%%%%%%%%%%%%%%%%%%%%%%%%%%%%%%%%%%%%%%%%
%%%%%%%%%%%%%%%%%%%%%%%%%%%%%%%%%%%%%%%%%%%%%%%%%%%%%%%%%%%%%%%%%%%%%%
%%%%%%%%%%%%%%%%%%%%%%%%%%%%%%%%%%%%%%%%%%%%%%%%%%%%%%%%%%%%%%%%%%%%%%
\subsection{Proof of Lemma~\ref{lem:xi_summary}}

As~\eqref{eq:xi} relates $\xi$ and $\chi$,
Item~\ref{item:xi:sum:1} in Lemma~\ref{lem:xi_summary} follows from Lemma~\ref{lem:chi:m=1}.
Item~\ref{item:xi:sum:2} and~\ref{item:xi:sum:4} are Corollary~\ref{cor:chi}.
We now prove Items~\ref{item:xi:sum:3} and~\ref{item:xi:sum:5}. Actually we provide a finer result.

For every $\theta\in (-1,1)$ and $k\geq 1$, by~\eqref{eq:sym:4} and~\eqref{eq:chi:6}, it holds
\begin{equation*}
\begin{split} 
	\frac{\xi_{2m}(\theta)}{(2m)!} 
	&= \chi_{2m+1}(\theta) 
	=\chi^{++}_{2m+1}(\theta) - \chi^{++}_{2m+1}(-\theta) + \chi^{+-}_{2m+1}(\theta) - \chi^{+-}_{2m+1}(-\theta)
	\\	
	& = \chi^{++}_{2m+1}(\theta) - \chi^{++}_{2m+1}(-\theta)  + \frac{1}{\sqrt{2\pi}} \frac{(1+\theta)^{2m+1}- (1-\theta)^{2m+1}}{(1-\theta^2)^{2m}}
\end{split}
\end{equation*}
and
\begin{equation*}
\begin{split} 
	- \frac{\xi_{2m-1}(\theta)}{(2m-1)!} 
	&= \chi_{2m}(\theta) 
	= \chi^{++}_{2m}(\theta) + \chi^{++}_{2m}(-\theta) + \chi^{+-}_{2m}(\theta) + \chi^{+-}_{2m}(-\theta)
	\\	
	& = \chi^{++}_{2m}(\theta) + \chi^{++}_{2m}(-\theta)  + \frac{1}{\sqrt{2\pi}} \frac{1+\theta}{(1-\theta)^{2m}} 
+  \frac{1}{\sqrt{2\pi}} \frac{1-\theta}{(1+\theta)^{2m}}.
\end{split}
\end{equation*}

In the proof of Lemma~\ref{lem:k:non_vanishing} we have already obtained a bound.
We find finer ones.
Let us recall the bound here: for every $\theta \in (0,1)$
\begin{equation*}
\begin{split} 
	\frac{\xi_{2m}(\theta)}{(2m)!} & = \chi_{2m+1}(\theta) 
	<
	- \frac{2\theta}{\sqrt{2\pi}} \Paren*{ \frac{1}{(1-\theta)^{2m}} -1}.
\end{split}
\end{equation*}
Similarly as in the proof of Lemma~\ref{lem:k:non_vanishing} we can obtain the desired bounds, 
by combining inequalities~\eqref{eq:chi:bounds} and Lemma~\ref{lem:theta:0}.
For every $\theta \in (0,1)$ we get
\begin{gather}
\begin{split} 
	\frac{\xi_{2m}(\theta)}{(2m)!} & 
	\geq - \frac{1+\sqrt{2\pi} \chi^{++}_{2m+1}(0)}{\sqrt{2\pi}} \Paren*{ \frac{1}{(1-\theta)^{2m}} - \frac{1}{(1+\theta)^{2m}}}
	\\
	& \qquad - \theta \frac{1-\sqrt{2\pi} \chi^{++}_{2m+1}(0)}{\sqrt{2\pi}} \Paren*{ \frac{1}{(1-\theta)^{2m}} + \frac{1}{(1+\theta)^{2m}}}
	\\
	& \geq - \frac{2}{\sqrt{2\pi}} \Paren*{ \frac{1}{(1-\theta)^{2m}} + \frac{1}{(1+\theta)^{2m}} } 
> - \frac{4}{\sqrt{2\pi}} \frac{1}{(1-\theta)^{2m}};
\end{split}
\\
\label{eq:xim:bounds:1}
\begin{split} 
	- \frac{\xi_{2m-1}(\theta)}{(2m-1)!} 
	&= \chi_{2m}(\theta) 
	\leq  \chi^{++}_{2m}(0) \Paren*{ 1+\theta +\frac{1}{(1-\theta)^{2(m-1)}}} 
	\\
	& \qquad + \frac{1}{\sqrt{2\pi}} \frac{1+\theta}{(1-\theta)^{2m-1}} +  \frac{1}{\sqrt{2\pi}} \frac{1-\theta}{(1+\theta)^{2m-1}},
\end{split}
\end{gather}
and
\begin{equation} 
    \label{eq:xim:bounds:2}
\begin{split} 
	- \frac{\xi_{2m-1}(\theta)}{(2m-1)!} 
	&	\geq  \chi^{++}_{2m}(0) \Paren*{ 1-\theta +\frac{1}{(1+\theta)^{2(m-1)}}} \\
	& \qquad + \frac{1}{\sqrt{2\pi}} \frac{1+\theta}{(1-\theta)^{2m-1}} +  \frac{1}{\sqrt{2\pi}} \frac{1-\theta}{(1+\theta)^{2m-1}}.
\end{split}
\end{equation}
We conclude that for all $\theta \in (-1,1)$,
\begin{equation*} 
	- \frac{4}{\sqrt{2\pi}}\frac1{(1-\abs{\theta})^{2m-1}}	<  \frac{\xi_{2m-1}(\theta)}{(2m-1)!} < - \frac{1+\abs{\theta}}{\sqrt{2\pi}}\frac1{(1-\abs{\theta})^{2m-1}}
\end{equation*}
as well as 
\begin{equation*}
\begin{split}
	& 	- \frac{4}{\sqrt{2\pi}} \frac{1}{(1-\abs{\theta})^{2m}}	
	< \frac{\xi_{2m}(\abs{\theta})}{(2m)!} 
	< - \frac{2\abs{\theta}}{\sqrt{2\pi}} \Paren*{ \frac{1}{(1-\abs{\theta})^{2m}} -1} \leq 0,
	\\
	& 0 \leq \frac{2\abs{\theta}}{\sqrt{2\pi}} \Paren*{ \frac{1}{(1-\abs{\theta})^{2m}} -1} 
	<  \frac{\xi_{2m}(-\abs{\theta})}{(2m)!}
	<  \frac{4}{\sqrt{2\pi}} \frac{1}{(1-\abs{\theta})^{2m}}.
\end{split}
\end{equation*}
These inequalities imply Items~\ref{item:xi:sum:3} and~\ref{item:xi:sum:5} of Lemma~\ref{lem:xi_summary} and complete its proof.

%%%%%%%%%%%%%%%%%%%%%%%%%%%%%%%%%%%%%%%%%%%%%%%%%%%%%%%%%%%%%%%%%%%%%%
%%%%%%%%%%%%%%%%%%%%%%%%%%%%%%%%%%%%%%%%%%%%%%%%%%%%%%%%%%%%%%%%%%%%%%
%%%%%%%%%%%%%%%%%%%%%%%%%%%%%%%%%%%%%%%%%%%%%%%%%%%%%%%%%%%%%%%%%%%%%%
\subsection{Proof of Proposition~\ref{prop:s_theta}}
\label{sec:amn}

We prove Proposition~\ref{prop:s_theta}
which completes the analysis of the asymptotic mixed normality in Theorem~\ref{thm:mixed_normality}.

The proof of Lemma~\ref{lem:xi_summary}, in particular~\eqref{eq:xim:bounds:1}-\eqref{eq:xim:bounds:2}, show that for all $\theta$,
\begin{equation*}
	- \frac{2}{\sqrt{2\pi}} \frac{1+\theta^2}{1-\theta^2} - \chi^{++}_2(0) (2+\abs{\theta}) 
	< \xi_1(\theta) 
	< - \frac{2}{\sqrt{2\pi}} \frac{1+\theta^2}{1-\theta^2} - \chi^{++}_2(0) (2-\abs{\theta}),
\end{equation*}
where $\chi^{++}_2(0) \approx 0.62$ by Lemma~\ref{lem:theta:0}.
This would already prove Proposition~\ref{prop:s_theta}, but the next lemma allow us to establish better bounds.

\begin{lemma} \label{lem:chi2}
Let 
\begin{equation*}
\mathcal{I}(\theta):= 
\int_0^\infty \int_0^\infty  \frac{(1-e^{-2 x y})}{(1-\theta^2 e^{-4 x y})} e^{-4 x y  - \frac{(x+y)^2}2} \vd x \vd y.
\end{equation*}
Then $\mathcal{I}$ is increasing on~$[0,1]$ and $\mathcal{I}(\theta) \in [\mathcal{I}(0), \mathcal{I}(1)]$ 
where
\begin{align*}
\mathcal{I}(0) %= \int_0^\infty \int_0^\infty  (1-e^{-2 x y}) e^{-4 x y  - \frac{(x+y)^2}2} dx dy
& = \sqrt{2\pi} (\chi^{++}_3(0) - \chi^{++}_4(0)) \approx  0.088,
\\
    \mathcal{I}(1)& =  \int_0^\infty \int_0^\infty  \frac{1}{(1+e^{-2 x y})} e^{-4 x y  - \frac{(x+y)^2}2} \vd x \vd y 
\in [\sqrt{2 \pi} \chi^{++}_3(0)/2 , \sqrt{2 \pi} \chi^{++}_4(0)],
\end{align*}
and
\begin{equation} \label{eq:chi2}
	\cs(\theta)^{-2}
	=-\xi_1(\theta)
	=\chi_2(\theta) 
	= \frac{2}{\sqrt{2\pi}} \left( \frac{1+\theta^2}{1-\theta^2} +  \sqrt{2\pi}\chi^{++}_2(0) - \theta^2 \mathcal{I}(\theta)\right).
\end{equation}
\end{lemma}

Recall that Lemma~\ref{lem:theta:0} shows that $\sqrt{2\pi} \chi^{++}_2(0)\approx 0.62$, $ \sqrt{2 \pi} \chi^{++}_3(0)/2 \approx 0.24 $, and $\sqrt{2 \pi} \chi^{++}_4(0) \approx 0.38.$

%%%
\begin{proof}
First note that $[-1,1] \ni \theta \mapsto \mathcal{I}(\theta)$ is an even function and increasing on~$[0,1]$, therefore the first part of the statement follows from~\eqref{eq:chi:3} and Lemma~\ref{lem:theta:0}.
We have already obtained the expressions~\eqref{eq:chi:6} and~\eqref{eq:chi:3bis} hence
\begin{multline*}
	\chi_2(\theta)  
	=
	\frac{1}{\sqrt{2\pi}} \frac{1+\theta^2}{1-\theta^2} 
	\\
	+  \int_0^\infty \int_0^\infty   \frac{ (1+ \theta) (1 -  \theta e^{- 2 x y} ) +  (1-\theta) (1 + \theta e^{- 2 x y} )}{ (1 -  \theta^2 e^{- 4 x y} )}  { e^{- 2 x y}}   e^{-\frac{(x+y)^2}2}  \vd y \vd x.
\end{multline*}
We rewrite the integrand in a suitable way
\begin{equation*}
\begin{split}
	\chi_2(\theta) 
	& = \frac{1}{\sqrt{2\pi}} \frac{1+\theta^2}{1-\theta^2} +  \frac2{\sqrt{2\pi}}\int_0^\infty \int_0^\infty  \frac{ (1- \theta^2 e^{- 2 x y} ) }{ (1 -  \theta^2 e^{- 4 x y} )}   { e^{- 2 x y}}   e^{-\frac{(x+y)^2}2} \vd y \vd x
	\\
	& = \frac{1}{\sqrt{2\pi}} \frac{1+\theta^2}{1-\theta^2} +  \frac2{\sqrt{2\pi}} \int_0^\infty \int_0^\infty  \left(1 - \frac{ \theta^2 e^{- 2 x y} (1-e^{- 2 x y}) }{ (1 -  \theta^2 e^{- 4 x y} )} \right)  { e^{- 2 x y}}  e^{-\frac{(x+y)^2}2} \vd y \vd x,
\end{split}
\end{equation*}
where we recognize the right-hand-side of~\eqref{eq:chi2}.
\end{proof}
%%%

%%%
\begin{remark} \label{rem:stheta}
We have obtained the following bounds for $\chi_2(\theta) = - \xi_1(\theta)=\cs(\theta)^{-2}$:
By \eqref{eq:chi2},
\begin{equation*}
	\frac{2}{\sqrt{2\pi}} \frac{1+\theta^2}{1-\theta^2} +  2 \chi^{++}_2(0)  - \frac{2}{\sqrt{2\pi}} \theta^2 \mathcal{I}(0)
	\leq  \chi_2(\theta) 
	\leq  \frac{2}{\sqrt{2\pi}} \frac{1+\theta^2}{1-\theta^2} + 2 \chi^{++}_2(0) - \frac{2}{\sqrt{2\pi}} \theta^2 \mathcal{I}(1).
\end{equation*}
Since $\mathcal{I}(0)= \sqrt{2\pi} (\chi^{++}_3(0) - \chi^{++}_4(0)) \approx  0.088$ 
and 
$\mathcal{I}(1)\leq \sqrt{2\pi} \chi^{++}_4(0)$, by studying the bounds above which are functions of $\theta$,
we can show that
\begin{equation*}	
    \frac{\cs(\theta)}{\sqrt{1-\theta^2}} \in [0.79, 0.88].
\end{equation*}
Numerically, we found that $\cs(\theta)$ is pretty close to its upper bound.
\end{remark}
%%%

%%%%%%%%%%%%%%%%%%%%%%%%%%%%%%%%%%%%%%%%%%%%%%%%%%%%%%%%%%%%%%%%%%%%%%
%%%%%%%%%%%%%%%%%%%%%%%%%%%%%%%%%%%%%%%%%%%%%%%%%%%%%%%%%%%%%%%%%%%%%%
%%%%%%%%%%%%%%%%%%%%%%%%%%%%%%%%%%%%%%%%%%%%%%%%%%%%%%%%%%%%%%%%%%%%%%

\subsection{Expansion around 0}
\label{sec:coeff_expansion}

We provide a polynomial expansion around $0$ of $\theta \mapsto \xi_m(\theta)$.
By~\eqref{eq:xi}, it suffices to provide it for $\theta \mapsto \chi_m(\theta)$.
Since we could not compute explicitly $\xi_m(\theta)$, this expansion allows us to compute it using their value at $\theta=0$ given in Lemma~\ref{lem:theta:0}.
Let us also recall that we could rewrite $\chi_m(\theta)$ in a form suitable for Monte Carlo simulation, see~\eqref{eq:MC}.

Let us define
\begin{equation*}
    \zeta_m(\theta):=\chi_m^{++}(\theta)+\chi_m^{+-}(\theta)-\chi_m^{--}(\theta)-\chi_m^{-+}(\theta). 
\end{equation*}
Using the symmetry relations \eqref{eq:sym:1}-\eqref{eq:sym:3}, $\theta\mapsto \zeta_{2m}(\theta)$ is odd, while $\theta\mapsto \zeta_{2m+1}(\theta)$ is
even. In particular, $\zeta_{2m}(0)=0$. 

Besides, 
\begin{equation*}
    \chi_{2m}(0)=\beta_{2m}
    \text{ and }
    \zeta_{2m+1}(0)=\beta_{2m+1}
\end{equation*}
with 
\begin{equation*}
    \beta_{m}:=\sqrt{\frac{2}{\pi}}\Paren*{\sqrt{2\pi}\chi^{++}_m(0)+(-1)^m}
\end{equation*}
and $\chi^{++}_m(0)$ provided in Lemma~\ref{lem:theta:0}.

\begin{lemma}
    For any $m\geq 1$, 
\begin{align}
    \label{eq:chi:deriv}
    \partial_\theta \chi_m(\theta)&=-(m-1)\chi_{m+1}(\theta)+\zeta_m(\theta), \\
    \label{eq:zeta:deriv}
    \partial_\theta \zeta_m(\theta)&=-(m-1)\zeta_{m+1}(\theta)+\chi_m(\theta). 
\end{align}
\end{lemma}

\begin{proof}
Differentiating with respect to $\theta$, 
\begin{multline*}
    \partial_\theta \chi^{\pm\pm}_m(\theta)
    =m\iint_{\RR^2}\ind{\pm x\geq 0,\pm y\geq 0}
   \partial_\theta k_\theta(x,y) k_\theta(x,y)^{m-1} p_\theta(1,x,y)\mu_\theta(x)\vd x\vd y
    \\
    +\iint_{\RR^2}\ind{\pm x\geq 0,\pm y\geq 0} k_\theta(x,y)^m \frac{\partial_\theta p_\theta(1,x,y)}{p_\theta(1,x,y)} p_\theta(1,x,y)\mu_\theta(x)\vd x\vd y
    \\
    +\iint_{\RR^2}\ind{\pm x\geq 0,\pm y\geq 0} \sgn(x) k_\theta(x,y)^m p_\theta(1,x,y)\vd x\vd y.
\end{multline*}

With \eqref{eq:derk} and since $k_\theta(x,y)=\partial_\theta p_\theta(1,x,y) /p_\theta(1,x,y)$,
\begin{multline*}
    \partial_\theta \chi^{\pm\pm}_m(\theta)
    =-m\iint_{\RR^2}\ind{\pm x\geq 0,\pm y\geq 0}
    \mu_\theta(x) k_\theta(x,y)^{m+1} p_\theta(1,x,y)\vd x\vd y
    \\
    +\iint_{\RR^2}\ind{\pm x\geq 0,\pm y\geq 0}
    \mu_\theta(x)k_\theta(x,y)^{m+1} p_\theta(1,x,y)\vd x\vd y
    \\
    +\iint_{\RR^2}\ind{\pm x\geq 0,\pm y\geq 0}\sgn(x)k_\theta(x,y)^m p_\theta(1,x,y)\vd x\vd y.
\end{multline*}
We then obtain
\begin{align*}
    \partial_\theta \chi^{++}_m(\theta)&=-(m-1) \chi^{++}_{m+1}(\theta)+\chi^{++}_m(\theta),\\
    \partial_\theta \chi^{+-}_m(\theta)&=-(m-1) \chi^{+-}_{m+1}(\theta)+\chi^{+-}_m(\theta),\\
    \partial_\theta \chi^{-+}_m(\theta)&=-(m-1) \chi^{-+}_{m+1}(\theta)-\chi^{-+}_m(\theta),\\
    \partial_\theta \chi^{--}_m(\theta)&=-(m-1) \chi^{--}_{m+1}(\theta)-\chi^{--}_m(\theta).
\end{align*}
With these four expressions, we obtain easily \eqref{eq:chi:deriv} and \eqref{eq:zeta:deriv}.
\end{proof}

\begin{corollary}
For $k\geq 1$ and $m\geq 1$, 
\begin{equation}
    \label{eq:chi:7}
    \partial^k_\theta\chi_m(\theta)
    =\sum_{\ell=0}^k (-1)^\ell b_\ell\binom{k}{\ell}
    \Paren[\Big]{
	\chi_{m+\ell}(\theta)a_\ell
	+
	\zeta_{m+\ell}(\theta)(1-a_\ell)
    }
\end{equation}
with 
\begin{equation*}
a_\ell
=
\begin{cases}
1&\text{ if }k-\ell\text{ is even},\\
0&\text{ if }k-\ell\text{ is odd}
\end{cases}
\text{ and }
b_\ell=
\begin{cases}
1&\text{ if }\ell=0,\\
(m-1)\dotsb (m-\ell-2)&\text{ if }\ell\geq 1.
\end{cases}.
\end{equation*}
\end{corollary}

\begin{proof}
We form the formal power series in $\RR[[\lambda,\mu]]$ for unknown $\lambda$ and $\mu$
by 
\begin{equation*}
X(\theta)=\sum_{m\geq 1}\lambda^m \chi_m(\theta)+\sum_{m\geq 1}\mu^m \zeta_m(\theta).
\end{equation*}
We then define two linear operators $S$ and $J$ on $\RR[[\lambda,\mu]]$ defined
by $S(\nu^{m+1})=-(m-1)\nu^m$ for $\nu=\lambda,\nu$
and $J(\lambda^m)=\mu^m$, $J(\mu^m)=\lambda^m$ so that $J$ exchanges $\mu$ and~$\lambda$
and~$J^2$ is idempotent. The operators~$S$ and~$J$ commute. 
With our choice of~$S$, 
\begin{equation*}
S(X(\theta))=-\sum_{m\geq 1}(m-1)\lambda^m\chi_{m+1}(\theta)-\sum_{m\geq 1}\mu^m \zeta_{m+1}(\theta),
\end{equation*}
so that we write \eqref{eq:chi:deriv}-\eqref{eq:zeta:deriv} as 
\begin{equation*}
\partial_\theta X(\theta)=(S+J) X(\theta).
\end{equation*}
Therefore, 
\begin{equation*}
    \partial^k_\theta X(\theta)=\sum_{\ell=0}^k \binom{k}{\ell}S^{\ell}J^{k-\ell} X(\theta).
\end{equation*}
We then deduce \eqref{eq:chi:7}.
\end{proof}

From \eqref{eq:chi:7},
\begin{align*}
    \partial^2_\theta \chi_m(\theta)&=(m-1)m\chi_{m+2}(\theta)-2(m-1)\zeta_{m+1}(\theta)+\chi_m(\theta),\\
    \partial^3_\theta \chi_m(\theta)&=-(m-1)m(m+1)\chi_{m+3}(\theta)+3(m-1)m\zeta_{m+2}(\theta)-3(m-1)\chi_{m+1}+\zeta_m(\theta).
\end{align*}

The Taylor expansion around $0$ leads to the following result.

\begin{corollary}
    For any $\theta$ with $\abs{\theta}<1$, $\chi_m(\theta)=\chi_m(0)+\sum_{k\geq 1}c_{m,k}\theta^k$
    where $c_{2m,1}=0$, $c_{2m+1,1}=-2m\chi_{2m+2}(0)$ and each of the $c_{m,k}$
    is a linear superposition of values of $\beta_\ell$ for $\ell=m,\dotsc,\ell+k$.
\end{corollary}

In particular, 
\begin{equation*}
    \chi_2(\theta)\approx 1.29+2.17\cdot\theta^2+\dotsb
    \text{ and }
    \chi_3(\theta)\approx -0.42\cdot\theta-2.93\cdot\theta^2+\dotsb.
\end{equation*}

%%%%%%%%%%%%%%%%%%%%%%%%%%%%%%%%%%%%%%%%%%%%%%%%%%%%%%%%%%%%%%%%%%%%%%
%%%%%%%%%%%%%%%%%%%%%%%%%%%%%%%%%%%%%%%%%%%%%%%%%%%%%%%%%%%%%%%%%%%%%%
%%%%%%%%%%%%%%%%%%%%%%%%%%%%%%%%%%%%%%%%%%%%%%%%%%%%%%%%%%%%%%%%%%%%%%
%%%%%%%%%%%%%%%%%%%%%%%%%%%%%%%%%%%%%%%%%%%%%%%%%%%%%%%%%%%%%%%%%%%%%%
\subsection{The limiting covariance}

\label{sec:Psi}
We now show how to compute, at least numerically, the covariance matrix~$\Psi$
appearing in Proposition~\ref{prop:score}.\ref{prop:score:i}. 
Recall that $\xi_{2i}(0)=0$  for any $i\geq 0$.

According to \cite[Theorem 1.2, p.~511]{jacod98} or \cite{mazzonetto19a}, 
\begin{multline*} 
    \Psi_{2i,2j} :=
    \int_{\RR} \EE\Prb*{\Paren*{\partial_\theta^{2i}k_0(B_0,B_1)\partial_\theta^{2j}k_0(B_0,B_1)}\given B_0=x}\vd x
    \\
+2\int_{\RR} \EE\Prb*{\sum_{\ell\geq 1}\partial_\theta^{2i}k_0(B_0,B_1)\partial_\theta^{2j}k_0(B_\ell,B_{\ell+1}) \given B_0=x}\vd x,
\end{multline*}
for a Brownian motion $B$.
Owing to \eqref{eq:derk} in Lemma~\ref{lem:k}, 
\begin{equation*} 
    \Psi_{2i,2j}
    =(2i)!(2j)!
    \int_{\RR} \EE\Prb*{\sum_{\ell\leq 0}k_0^{2i+1}(B_0,B_1)k^{2j+1}_0(B_\ell,B_{\ell+1})
\given B_0=x}\vd x.
\end{equation*}
We decompose this expression as 
\begin{equation*}
    \Psi_{2i,2j} = (2i)!(2j)!\Paren*{\uPsi^I_{2i,2j} + \uPsi^{II}_{2i,2j} + \uPsi^{III}_{2i,2j}}
\end{equation*}
with 
\begin{align*}
    \uPsi^I_{2i,2j}
    =&
    \int_{\RR} \EE\Prb*{k_0^{2i+1}(B_0,B_1)k_0^{2j+1}(B_0,B_1) \given B_0=x}\vd x,
\\
\uPsi^{II}_{2i,2j}
    =&
    \int_{\RR} \EE\Prb*{k_0^{2i+1}(B_0,B_1)k^{2j+1}_0(B_1,B_2) \given B_0=x}\vd x
    \\
    &+\int_{\RR} \EE\Prb*{k_0^{2j+1}(B_0,B_1)k^{2i+1}_0(B_1,B_2) \given B_0=x}\vd x,
\\
\uPsi^{III}_{2i,2j}
    =&
    \int_{\RR} \EE\Prb*{\sum_{\ell\geq 2}k_0^{2i+1}(B_0,B_1)k^{2j+1}_0(B_\ell,B_{\ell+1}) \given B_0=x}\vd x
    \\
     &+\int_{\RR} \EE\Prb*{\sum_{\ell\geq 2}k_0^{2j+1}(B_0,B_1)k^{2i+1}_0(B_\ell,B_{\ell+1}) \given B_0=x}\vd x.
\end{align*}

\begin{lemma}
    \label{lemma:Psi:m0}
It holds that
    \begin{equation}
	\label{eq:Psi:I}
	\uPsi_{2i,2j}^I = \chi_{2(i+j)+2}(0) = \frac{2}{\sqrt{2\pi}} \left( 1+ \sqrt{2\pi}\chi^{++}_{2(i+j+1)}(0)\right) = \frac{-\xi_{2(i+j)+1}(0)}{(2(i+j)+1)!}.
    \end{equation}
    In particular, $\Xi_{0,0}(0)=\sqrt{\uPsi^I_{0,0}}=\sqrt{-\xi_1(0)}$.
\end{lemma}
\begin{proof}
Eq.~\eqref{eq:Psi:I} follows from the definition of $\xi_m$ in~\eqref{eq:xi:def} and \eqref{eq:xi}.
When $i=j=0$, by the Markov property and~\eqref{eq:average:0}, we have $\uPsi^{II}_{0,0}=\uPsi^{III}_{0,0}=0$.
\end{proof}

By~\eqref{eq:chi:decomp} and Lemma~\ref{lem:chimp}, $\chi_{2(i+j)+2}(0)= \frac{2}{\sqrt{2\pi}} (1+\sqrt{2\pi}\chi_{2(i+j)+2}^{++}(0))$, and an approximation of $\chi^{++}(0)$ is given in~Lemma~\ref{lem:theta:0}.

For $m\geq 1$, $m$ odd, we set 
\begin{align*}
    \kappa_m(x,y) &:= \exp(-2m(xy)_+)\sgn(y)=k_0(x,y)^m
    \\
    \text{ and }
    \widehat{\kappa}_m(x,y) &:= \exp(-2m(xy)_+)\sgn(x)=\kappa_m(y,x).
\end{align*}
We also define
\begin{align*}
    K_m(x) &:= \int \kappa_m(x)\exp\Paren*{\frac{-(x-y)^2}{2}}\frac{\vd y}{\sqrt{2\pi}}
    =\EE\Prb{\kappa_m(x,B_1)\given B_0=x}
    \\
    \text{ and }
    \widehat{K}_m(x) &:= \int \widehat{\kappa}_m(x)\exp\Paren*{\frac{-(x-y)^2}{2}}\frac{\vd y}{\sqrt{2\pi}}
=\EE\Prb{\widehat{\kappa}_m(x,B_1)\given B_0=x}
\end{align*}
for a Brownian motion $B$.
With some straightforward computations, 
\begin{align*}
    K_m(x)&=\sgn(x)\exp(2x^2m(m-1))\Phi(-(2m-1)\abs{x})
-\sgn(x)\Phi(-\abs{x})
\\
\text{ and }
    \widehat{K}_m(x)&=\sgn(x)\exp(2x^2m(m-1))\Phi(-(2m-1)\abs{x})
+\sgn(x)\Phi(-\abs{x}).
\end{align*}

\begin{notation}[(Scaled) complementary error function]
    The \emph{complementary error function}
    and the \emph{scaled complementary error function} are
    \begin{equation*}
	\erfc(x):=\frac{2}{\sqrt{\pi}}\int_x^{+\infty} e^{-y^2}\vd y
	=2\overline{\Phi}(x\sqrt{2})
	\text{ and }
	\erfcx(x):=\exp\Paren*{x^2}\erfc(x)
    \end{equation*}
    for $x\geq 0$.
\end{notation}
The functions $x\in[0,1]\mapsto\erfc(x)$ and $x\in[0,1]\mapsto\erfcx(x)$ 
take their values in~$[0,1]$.
Thanks to the Mill's ratio, for $x$ large, $\erfcx(x)\sim 1/x\sqrt{\pi}$. 

We also introduce 
\begin{equation*}
    b_m:=\frac{2m-1}{\sqrt{2}}\text{ so that }2m(m-1)=b_m^2-\frac{1}{2}.
\end{equation*}
As $\Phi(-\abs{x})=\overline{\Phi}(\abs{x})$, we rewrite $K_m$ and $\widehat{K}_m$
as
\begin{align*}
    K_m(x)&=\frac{1}{2}\sgn(x)e^{-x^2/2}
    \Paren*{\erfcx(b_m\abs{x})-\erfcx(\abs{x}/\sqrt{2})}
    \\
    \text{ and }
    \widehat{K}_m(x)&=\frac{1}{2}\sgn(x)e^{-x^2/2}
    \Paren*{\erfcx(b_m\abs{x})+\erfcx(\abs{x}/\sqrt{2})}.
\end{align*}

%%%%
\begin{lemma}
    For $i,j\geq 0$, 
    \begin{multline*}
	\uPsi^{III}_{2i,2j}
	=
	\sum_{\ell\geq 2} 
	\frac{1}{\sqrt{2\pi}\sqrt{\ell-1}}
    \int_{-\infty}^{+\infty}\int_{-\infty}^{+\infty}
    \widehat{K}_{2i}(y)K_{2j}(z)\exp\Paren*{\frac{-(z-y)^2}{2(\ell-1)}}\vd y\vd z
    \\
    +\sum_{\ell\geq 2}
    \frac{1}{\sqrt{2\pi}\sqrt{\ell-1}}
    \int_{-\infty}^{+\infty}\int_{-\infty}^{+\infty}
    \widehat{K}_{2j}(y)K_{2i}(z)\exp\Paren*{\frac{-(z-y)^2}{2(\ell-1)}}\vd y\vd z.
\end{multline*}
\end{lemma}
%%%%

\begin{proof}
We define for $m,n\geq0$,
\begin{equation}
    \label{eq:Psi:A}
    A(\ell,m,n) := \int_{-\infty}^{+\infty} \EE\Prb{\kappa_m(B_0,B_1)\kappa_n(B_\ell,B_{\ell+1})\given B_0=x}\vd x
\end{equation}
so that 
\begin{equation*}
    \uPsi^{III}_{2i,2j}=\sum_{\ell\geq 2} 2A(\ell,2i+1,2j+1).
\end{equation*}
We rewrite $A(\ell,m,n)$ using its integral form:
\begin{multline}
    A(\ell,m,n)
    =
    \frac{1}{(2\pi)^{3/2}\sqrt{\ell-1}}\int_{-\infty}^{+\infty}\int_{-\infty}^{+\infty}\int_{-\infty}^{+\infty}\int_{-\infty}^{+\infty} \kappa_m(x,y)\kappa_n(z,u)
    \\
    \times
    \exp\Paren*{\frac{-(x-y)^2}{2}}
    \exp\Paren*{\frac{-(z-y)^2}{2(\ell-1)}}
    \exp\Paren*{\frac{-(z-u)^2}{2}}
    \vd x\vd y\vd z\vd u.
\end{multline}

We compute first $A(\ell,m,n)$ for $\ell\geq 2$.
Since $\widehat{\kappa}_m(y,x)=\kappa_m(x,y)$, inverting~$x$ and~$y$,
and using the definitions of $K_m$ and $\widehat{K}_m$ leads to 
\begin{equation*}
    A(\ell,m,n)=\frac{1}{\sqrt{2\pi}\sqrt{\ell-1}}
    \int_{-\infty}^{+\infty}\int_{-\infty}^{+\infty}
    \widehat{K}_m(y)K_n(z)\exp\Paren*{\frac{-(z-y)^2}{2(\ell-1)}}\vd y\vd z.
\end{equation*}
This gives the result.
\end{proof}

\begin{remark}
With a change of variable $(y,z)\to (y/\sqrt{j-1},z/\sqrt{j-1})$, 
for $A$ introduced in \eqref{eq:Psi:A}, 
\begin{multline*}
    A(\ell,m,n)=\sqrt{2\pi}\sqrt{\ell-1}
    \int_{-\infty}^{+\infty}\int_{-\infty}^{+\infty}
    \widehat{K}_m(y\sqrt{\ell-1})K_n(z\sqrt{\ell-1})
    \exp\Paren{yz}
    \\
    \times 
    \exp\Paren*{\frac{-z^2}{2}}
    \exp\Paren*{\frac{-y^2}{2}}\frac{\vd y}{\sqrt{2\pi}}\frac{\vd z}{\sqrt{2\pi}}.
\end{multline*}
For two independent Gaussian random variables $G,G'\sim\cN(0,1)$, we may write 
\begin{equation*}
    A(\ell,m,n) = \sqrt{2\pi}\sqrt{\ell-1}\EE\Prb{K_m\Paren{G'\sqrt{\ell-1}}\widehat{K}_n\Paren{G\sqrt{\ell-1}}
    \exp(G\cdot G')},
\end{equation*}
where $G,G'\sim\cN(0,1)$ are independent. However, the variance is high so that Monte Carlo 
techniques are unstable here.
\end{remark}
%%

%%%%
\begin{lemma}
    For $i,j\geq 1$, 
\begin{equation*}
    \uPsi^{II}_{2i,2j}=
    \int_0^{+\infty} e^{-x^2}\erfcx(b_{2i+1}x)\erfc(b_{2j+1}x)\vd x
-\frac{2}{\sqrt{\pi}}\Paren*{\sqrt{2}-1}.
\end{equation*}
\end{lemma}
%%%%

\begin{proof}
As above, 
\begin{multline*}
    A(1,m,n)=\int_{-\infty}^{+\infty} \widehat{K}_m(x)K_n(x)\vd x
    \\
    =
    \frac{1}{2}\int_0^{+\infty} e^{-x^2}\erfcx(b_mx)\erfc(b_n x)\vd x
    -\frac{1}{2}\int_0^{+\infty} e^{-x^2}\erfcx(x/\sqrt{2})^2\vd x.
\end{multline*}
The quantity $\uPsi_{2i,2j}^{III}$ is then 
\begin{equation*}
    \uPsi_{2i,2j}^{II} = A(1,2i+1,2j+1) + A (1,2j+1,2i+1).
\end{equation*}
According to  \cite[4.7.2]{Ng1},
\begin{equation*}
    \int_0^{+\infty} e^{-x^2}\erfcx\Paren*{\frac{x}{\sqrt{2}}}^2\vd x
=
\int_0^{+\infty} \erfc\Paren*{\frac{x}{\sqrt{2}}}^2\vd x
=\frac{2}{\sqrt{\pi}}\Paren*{\sqrt{2}-1}
\approx 0.467.
\end{equation*}
No closed-form expression for $\int_{-\infty}^{+\infty} e^{-x^2}\erfcx(b_mx)\erfc(b_n x)\vd x$
seems to exists. However, this integral is easy to compute numerically
as there exists various implementation of $\erfcx$.
This could be done through a quadrature method to compute the integral
or using a Monte Carlo method as 
\begin{multline*}
\int_0^{+\infty} e^{-x^2}\erfcx(b_mx)\erfc(b_nx)\vd x
\\
=
\frac{\sqrt{\pi}}{\sqrt{2m-1}}\EE\Prb*{\erfcx\Paren*{\frac{\sqrt{2m-1}}{\sqrt{2}}\abs{G}}\erfcx\Paren*{\frac{\sqrt{2n-1}}{\sqrt{2}}\abs{G}}}.
\end{multline*}
This concludes the proof.
\end{proof}

\paragraph*{Numerical computations.} 
We give some the first values of $\Psi$
obtained by numerical computations:
\begin{equation*}
    \Set{\Psi_{2i,2j}}_{i,j=0,1,2}
    \approx
    \begin{bmatrix}
	\phantom{00}1.295 & \phantom{00}1.422 & \phantom{0}12.576 \\
	\phantom{00}1.422 & \phantom{00}1.891 & \phantom{0}18.135 \\
	\phantom{0}12.576 & \phantom{0}18.135 & 181.421 
    \end{bmatrix}.
\end{equation*}

Let us conclude by an observation: for $m\geq 2$, $m$ odd, 
we found the empirical rule of thumbs that $\uPsi_{m,m}/\uPsi^I_{m,m}$
behaves, for moderate values of $m$, as $C_m/\sqrt{2m}$
for a pre-factor $C_m$ varying slowly. In addition, 
$\uPsi^I_{m,m}$ 
converges to $\frac{2}{\sqrt{2\pi}}$
as $m\to\infty$ (see Lemma~\ref{lem:theta:0} and Lemma~\ref{lemma:Psi:m0})
while $\uPsi^{II}_{m,m}$ converges to $-2(\sqrt{2}-1)/\sqrt{\pi}$.

%%%%%%%%%%%%%%%%%%%%%%%%%%%%%%%%%%%%%%%%%%%%%%%%%%%%%%%%%%%%%%%%%%%%%%
%%%%%%%%%%%%%%%%%%%%%%%%%%%%%%%%%%%%%%%%%%%%%%%%%%%%%%%%%%%%%%%%%%%%%%
%%%%%%%%%%%%%%%%%%%%%%%%%%%%%%%%%%%%%%%%%%%%%%%%%%%%%%%%%%%%%%%%%%%%%%
%%%%%%%%%%%%%%%%%%%%%%%%%%%%%%%%%%%%%%%%%%%%%%%%%%%%%%%%%%%%%%%%%%%%%%
\section{Proofs of the asymptotic results}

\label{sec:proofs}

The goal of this section is to prove Proposition~\ref{prop:score}.
Proposition~\ref{prop:pn_dn} is proved in~Section~\ref{sec:score} as a
corollary of Proposition~\ref{prop:score}.  Theorem~\ref{thm:asym} is a direct
application of~\cite[Theorem 3]{lm22}.  The convergence in
Theorem~\ref{thm:mixed_normality} follows from combining the expansion in
Theorem~\ref{thm:asym} with Proposition~\ref{prop:pn_dn}.

%%%%%%%%%%%%%%%%%%%%%%%%%%%%%%%%%%%%%%%%%%%%%%%%%%%%%%%%%%%%%%%%%%%%%%
%%%%%%%%%%%%%%%%%%%%%%%%%%%%%%%%%%%%%%%%%%%%%%%%%%%%%%%%%%%%%%%%%%%%%%
%\subsection{Integrability conditions}

Proposition~\ref{prop:score} is a consequence of the results stated in
Remark~\ref{rem:score}, in particular of the convergence~\eqref{eq:rate:conv}
which is a simple application of~\cite{mazzonetto19a}.  The next
Lemma~\ref{lem:int:k} is an integrability condition which ensures that the
assumptions of the main result in~\cite{mazzonetto19a} are satisfied.
Although~\eqref{eq:rate:conv} concerns derivatives of $k_\theta$, the lemma
concerns its powers. This is justified by Lemma~\ref{lem:k} in
Section~\ref{sec:score_coefficients} where it is shown that 
$\partial_\theta^m k_\theta(x,y)$ is proportional to $k_\theta(x,y)^m$ with a multiplicative
factor depending only on $m\in \NN$. 

%%%% INTEGRABILITY

\begin{lemma} 
    \label{lem:int:k}
	 For $m=1,2,3,\dotsc,$ and any $\gamma>3$, there exists $a>0$ 
	 such that for every $\theta\in (-1,1)$ there exists 
    a measurable, bounded function ${h}_\theta:\RR\to\RR$ 
		which satisfies
    \begin{equation*}
    \int_{\RR} {h}_\theta(x)\abs{x}^\gamma\vd x<+\infty
    \text{ and }
	\abs{k_\theta^m(x,y)}\leq {h}_\theta (x)e^{a\abs{y-x}},\text{ for any }x,y\in\RR.
    \end{equation*}
\end{lemma}

\begin{proof}
Fix $a>0$.  Observe that 
\begin{equation*}
    \abs{ k_\theta(x,y) }
	\leq  
	\frac{1}{(1-\abs{\theta})} \ind{\{x y \leq 0\}} 
	+  \frac1{(1-\abs{\theta})} e^{-2 x y} \ind{\{x y > 0\}} =: h_\theta(x,y) \leq \frac{1}{(1-\abs{\theta})}.
\end{equation*}
Now, 
$\ind{\Set{x y \leq 0}} \leq e^{-a \abs{x}} e^{a \abs{y-x}}$. 
Let us find a bound for
$e^{- 2 x y} e^{- a\abs{y - x}} \ind{\Set{x y > 0}}$.
Let us first assume $x y>0$ and $|y| \geq |x|$ then
\begin{equation*}
	e^{- 2 x y} e^{- a|y - x|} 
	= e^{- 2 |x| |y|} e^{- a|y| + a |x|} 
	\leq  e^{- 2 x^2}
    \end{equation*}
and if $x y>0$ and $\abs{y} < \abs{x}$ then
\begin{equation*}
	e^{- 2 x y} e^{- a|y - x|} 
	= e^{- 2 \abs{x} \abs{y} } e^{- a\abs{x} + a \abs{y}}
	\leq  e^{- 2 x^2} \ind{\{\abs{x}\leq a /2\}} + e^{- a\abs{x}} \ind{\{\abs{x} > a /2\}}.
    \end{equation*}
We conclude that for any $x,y\in\RR$, 
\begin{equation*}
    \abs{k_\theta(x,y)}^m \leq  (1-\abs{\theta})^{-(m-1)}  h_\theta(x,y) \leq (1-\abs{\theta})^{-m} \left( e^{-2 x^2} + e^{-a \abs{x}}\right) e^{a\abs{y-x}}
\end{equation*}
with $\int_{\RR} \Paren[\big]{ e^{-2 x^2} + e^{-a |x|}}\abs{x}^\gamma\vd x<+\infty$. 
\end{proof}
%%%%%

%%%%%%%%%%%%%%%%%%%%%%%%%%%%%%%%%%%%%%%%%%%%%%%%%%%%%%%%%%%%%%%%%%%%%%
%%%%%%%%%%%%%%%%%%%%%%%%%%%%%%%%%%%%%%%%%%%%%%%%%%%%%%%%%%%%%%%%%%%%%%
%\subsection{Proof of Proposition~\ref{prop:score}}

%%%%%%%%%%%%%%%%%%%%%%%%%%%%%%%%%%%%%%%%%%%%%%%%%%%%%%%%%%%%%%%%%%%%%%
\subsubsection{Item~\ref{prop:score:ii}: Case \texorpdfstring{$\theta\neq0$}{θ≠0}}

Item~\ref{prop:score:ii} of Proposition~\ref{prop:score} is an 
immediate consequence of the result stated in Remark~\ref{rem:score}, 
Remark~\ref{rem:rate:conv}, and 
Lemma~\ref{lem:xi_summary} which ensures that
when $\theta\in(-1,1)$,
$\xi_0(\theta)=0$,
and
when in addition $\theta\neq0$,
$\xi_m(\theta)\neq 0$  as soon as $m\geq 1$
and $\xi_{2m}(0)=0$.

Using the fact that given two sequences, one converging in probability and the other converges stably, 
joint stable convergence holds (see~\cite[Theorem~1]{zbMATH03586244})), we get the joint
stable convergence of the vector $(\coef{S}_k(\theta))_{k=0,\dotsc,m}$ for any $m\geq 0$.

%%%%%%%%%%%%%%%%%%%%%%%%%%%%%%%%%%%%%%%%%%%%%%%%%%%%%%%%%%%%%%%%%%%%%%
\subsubsection{Item~\ref{prop:score:i}: Case \texorpdfstring{$\theta=0$}{θ=0}}

\label{sec:Cov}

From the above arguments 
we have only to focus on the convergence of
\begin{equation*}
\coef{S}^{\zig}(n,0):=(n^{1/4}\coef{S}_0(n,0),n^{1/4}\coef{S}_2(n,0),\dotsc,n^{1/4}\coef{S}_{2m}(n,0))
\end{equation*}
since $\coef{S}_{2m+1}(n,0)$ converges in probability to $\xi_{2m+1}(\theta)L_T$
for each $m\geq 0$.

The result is a direct consequence of~\cite[Theorem 1.2, p.~511]{jacod98} (which can be applied since Lemma~\ref{lem:int:k} holds).

%%%%%%%%%%%%%%%%%%%%%%%%%%%%%%%%%%%%%%%%%%%%%%%%%%%%%%%%%%%%%%%%%%%%%%
%%%%%%%%%%%%%%%%%%%%%%%%%%%%%%%%%%%%%%%%%%%%%%%%%%%%%%%%%%%%%%%%%%%%%%
%%%%%%%%%%%%%%%%%%%%%%%%%%%%%%%%%%%%%%%%%%%%%%%%%%%%%%%%%%%%%%%%%%%%%%

%%%%%%%%%%%%%%%%%%%%%%%%%%%%%%%%%%%%%%%%%%%%%%%%%%%%%%%%%%%%%%%%%%%%%%
%%%%%%%%%%%%%%%%%%%%%%%%%%%%%%%%%%%%%%%%%%%%%%%%%%%%%%%%%%%%%%%%%%%%%%
% BIBLIOGRAPHY
%%%%%%%%%%%%%%%%%%%%%%%%%%%%%%%%%%%%%%%%%%%%%%%%%%%%%%%%%%%%%%%%%%%%%%
%%%%%%%%%%%%%%%%%%%%%%%%%%%%%%%%%%%%%%%%%%%%%%%%%%%%%%%%%%%%%%%%%%%%%%

%\bibliographystyle{plain}
%\bibliography{biblio-rate-sbm.bib}

\begin{thebibliography}{10}

\bibitem{abramowitz}
Milton Abramowitz and Irene~A. Stegun.
\newblock {\em Handbook of mathematical functions with formulas, graphs, and
  mathematical tables}, volume~55 of {\em National Bureau of Standards Applied
  Mathematics Series}.
\newblock 1964.

\bibitem{zbMATH03586244}
D.~J. Aldous and G.~K. Eagleson.
\newblock On mixing and stability of limit theorems.
\newblock {\em Ann. Probab.}, 6:325--331, 1978.

\bibitem{araya}
H\'{e}ctor Araya, Meryem Slaoui, and Soledad Torres.
\newblock Bayesian inference for fractional oscillating {B}rownian motion.
\newblock {\em Comput. Statist.}, 37(2):887--907, 2022.

\bibitem{Bai2021}
Yizhou Bai, Yongjin Wang, Haoyan Zhang, and Xiaoyang Zhuo.
\newblock {B}ayesian estimation of the skew {O}rnstein-{U}hlenbeck process.
\newblock {\em Computational Economics}, 60(2):479--527, July 2021.

\bibitem{barahona}
Manuel Barahona, Laura Rifo, Maritza Sep\'{u}lveda, and Soledad Torres.
\newblock A simulation-based study on {B}ayesian estimators for the skew
  {B}rownian motion.
\newblock {\em Entropy}, 18(7):Paper No. 241, 14, 2016.

\bibitem{bardou}
Olivier Bardou and Miguel Martinez.
\newblock Statistical estimation for reflected skew processes.
\newblock {\em Stat. Inference Stoch. Process.}, 13(3):231--248, 2010.

\bibitem{bounebache_2014}
Said~Karim Bounebache and Lorenzo Zambotti.
\newblock A skew stochastic heat equation.
\newblock {\em J. Theor. Probab.}, 27(1):168--201, 2014.

\bibitem{cantrell99a}
R.S. Cantrell and C.~Cosner.
\newblock Diffusion models for population dynamics incorporating individual
  behavior at boundaries: Applications to refuge design.
\newblock {\em Theoretical Population Biology}, 55(2):189--207, 1999.

\bibitem{decamps}
Marc Decamps, Marc Goovaerts, and Wim Schoutens.
\newblock Self exciting threshold interest rates models.
\newblock {\em Int. J. Theor. Appl. Finance}, 9(7):1093--1122, 2006.

\bibitem{gairat17a}
Alexander Gairat and Vadim Shcherbakov.
\newblock Density of skew {B}rownian motion and its functionals with
  application in finance.
\newblock {\em Math. Finance}, 27(4):1069--1088, 2017.

\bibitem{harrison}
J.~M. Harrison and L.~A. Shepp.
\newblock On skew {B}rownian motion.
\newblock {\em Ann. Probab.}, 9(2):309--313, 1981.

\bibitem{ih}
I.A. Ibragimov and R.Z. Has'minskii.
\newblock {\em Statistical Estimation Asymptotic Theory}.
\newblock Springer, 1981.

\bibitem{ILM_22}
Andrey Itkin, Alexander Lipton, and Dmitry Muravey.
\newblock Multilayer heat equations and their solutions via oscillating
  integral transforms.
\newblock {\em Physica A}, 601:34, 2022.
\newblock Id/No 127544.

\bibitem{jacod98}
J.~Jacod.
\newblock Rates of convergence to the local time of a diffusion.
\newblock {\em Ann. Inst. H. Poincar\'e Probab. Statist.}, 34(4):505--544,
  1998.

\bibitem{K}
Yu.A. Kutoyants.
\newblock {\em Parameter Estimation for Stochastic Processes}.
\newblock Heldermann, 1984.

\bibitem{lejay_sbm}
A.~Lejay.
\newblock On the constructions of the skew {B}rownian motion.
\newblock {\em Probab. Surv.}, 3:413--466, 2006.

\bibitem{pigato}
A.~Lejay and P.~Pigato.
\newblock Statistical estimation of the oscillating brownian motion.
\newblock {\em Bernoulli}, 24(4B):3568--3602, 2018.

\bibitem{lejay2018}
Antoine Lejay.
\newblock Estimation of the bias parameter of the skew random walk and
  application to the skew {B}rownian motion.
\newblock {\em Stat. Inference Stoch. Process.}, 21(3):539--551, 2018.

\bibitem{lejay18}
Antoine Lejay.
\newblock Estimation of the bias parameter of the skew random walk and
  application to the skew {B}rownian motion.
\newblock {\em Stat. Inference Stoch. Process.}, 21(3):539--551, 2018.

\bibitem{lm22}
Antoine Lejay and Sara Mazzonetto.
\newblock Beyond the delta method, 2022.
\newblock Preprint \url{https://hal.inria.fr/hal-03738371v1}.

\bibitem{lejay2014}
Antoine Lejay, Ernesto Mordecki, and Soledad Torres.
\newblock Is a {Brownian} motion skew?
\newblock {\em Scand. J. Stat.}, 41(2):346--364, 2014.

\bibitem{lejay2019}
Antoine Lejay, Ernesto Mordecki, and Soledad Torres.
\newblock Two consistent estimators for the skew brownian motion.
\newblock {\em ESAIM: Probability and Statistics}, 23:567--583, 2019.

\bibitem{lejay12}
Antoine Lejay and G\'{e}raldine Pichot.
\newblock Simulating diffusion processes in discontinuous media: a numerical
  scheme with constant time steps.
\newblock {\em J. Comput. Phys.}, 231(21):7299--7314, 2012.

\bibitem{lejay_pigato2018}
Antoine Lejay and Paolo Pigato.
\newblock Statistical estimation of the oscillating {B}rownian motion.
\newblock {\em Bernoulli}, 24(4B):3568--3602, 2018.

\bibitem{lejay_pigato2019}
Antoine Lejay and Paolo Pigato.
\newblock A threshold model for local volatility: evidence of leverage and mean
  reversion effects on historical data.
\newblock {\em Int. J. Theor. Appl. Finance}, 22(4):1950017, 24, 2019.

\bibitem{lejay_pigato2020}
Antoine Lejay and Paolo Pigato.
\newblock Maximum likelihood drift estimation for a threshold diffusion.
\newblock {\em Scand. J. Stat.}, 47(3):609--637, 2020.

\bibitem{mazzonetto19a}
Sara Mazzonetto.
\newblock Rates of convergence to the local time of oscillating and skew
  brownian motions, 2019.

\bibitem{mota14a}
Pedro~P. Mota and Manuel~L. Esqu{\'{\i}}vel.
\newblock On a continuous time stock price model with regime switching, delay,
  and threshold.
\newblock {\em Quant. Finance}, 14(8):1479--1488, 2014.

\bibitem{Ng1}
Edward~W. Ng and Murray Geller.
\newblock A table of integrals of the error functions.
\newblock {\em J. Res. Nat. Bur. Standards Sect. B}, 73B:1--20, 1969.

\bibitem{ovaskainen03a}
Otso Ovaskainen and Stephen~J. Cornell.
\newblock Biased movement at a boundary and conditional occupancy times for
  diffusion processes.
\newblock {\em J. Appl. Probab.}, 40(3):557--580, 2003.

\bibitem{ramirez2}
Jorge~M. Ramirez, Enrique~A. Thomann, Edward~C. Waymire, Roy Haggerty, and
  Brian Wood.
\newblock A generalized {T}aylor-{A}ris formula and skew diffusion.
\newblock {\em Multiscale Model. Simul.}, 5(3):786--801, 2006.

\bibitem{robert_2022}
Christian~Y. Robert.
\newblock How large is the jump discontinuity in the diffusion coefficient of a
  time-homogeneous diffusion?
\newblock {\em Econometric Theory}, page 1–33, 2022.

\bibitem{su_chan2015}
Fei Su and Kung-Sik Chan.
\newblock Quasi-likelihood estimation of a threshold diffusion process.
\newblock {\em J. Econometrics}, 189(2):473--484, 2015.

\bibitem{su_chan2017}
Fei Su and Kung-Sik Chan.
\newblock Testing for threshold diffusion.
\newblock {\em J. Bus. Econom. Statist.}, 35(2):218--227, 2017.

\bibitem{walsh}
J.~B. Walsh.
\newblock A diffusion with discontinuous local time.
\newblock In {\em Temps locaux}, volume 52-53, pages 37--45. Soci{\'e}t{\'e}
  Math{\'e}matique de France, 1978.

\bibitem{zhang}
M.~Zhang.
\newblock Calculation of diffusive shock acceleration of charged particles by
  skew {B}rownian motion.
\newblock {\em The Astrophysical Journal}, 541:428--435, 2000.

\end{thebibliography}

\end{document}